\documentclass[english]{article}
\usepackage[T1]{fontenc}
\usepackage[latin9]{inputenc}
\usepackage{color}
\usepackage{verbatim}
\usepackage{amsthm}
\usepackage{amsmath}
\usepackage{amssymb}
\usepackage{esint}
\makeatletter
\numberwithin{equation}{section}
\numberwithin{figure}{section}

\theoremstyle{plain}
\newtheorem{thm}{Theorem}
  \theoremstyle{definition}
\newtheorem{defn}[thm]{Definition}
  \theoremstyle{plain}
\newtheorem{lem}[thm]{Lemma}
  \theoremstyle{plain}
\newtheorem{prop}[thm]{Proposition}
  \theoremstyle{remark}
\newtheorem{rem}[thm]{Remark}
 \theoremstyle{definition}
\newtheorem{example}[thm]{Example}

\usepackage[all]{xy}
\usepackage{amscd}
\makeatletter 
\newcommand{\xyR}[1]{ \makeatletter
\xydef@\xymatrixrowsep@{#1} \makeatother} 
\newcommand{\xyC}[1]{ \makeatletter
\xydef@\xymatrixcolsep@{#1} \makeatother} 
\entrymodifiers={++[ ][F]} 
\newcommand{\h}[1]{\hspace{#1 truein}}  

\newcommand{\freccia}{\longrightarrow} 

\newcommand{\eps}{\varepsilon} 
\renewcommand{\phi}{\varphi} 
\newcommand{\field}[1]{\mathbb{#1}}
\newcommand{\R}{\field{R}}                        
\newcommand{\ER}{{^\bullet\R}}                    
\newcommand{\N}{\field{N}}                        
\newcommand{\hyperR}{{^*\R}}                      
\newcommand{\diff}[1]{\,{\rm d}#1}

\newcommand{\st}[1]{{^\circ #1}} 
\newcommand{\Cc}{\mathcal{C}} 
\newcommand{\D}{\mathcal{D}} 
\newcommand{\ext}[1]{{}^\bullet #1} 

\DeclareMathOperator*{\ManInfty}{{\bf Man}}

\newcommand{\CInfty}{\boldsymbol{\mathcal{C}}^{\infty}} 
\newcommand{\ECInfty}{\ext{\CInfty}} 

\newcommand{\ptind}{\displaystyle \mathop {\ldots\ldots\,}} 
\newcommand{\pti}{:\;\;\;} 
\newcommand{\then}{\quad \Longrightarrow \quad}

\makeatother

\renewcommand{\qed}{\hspace*{\fill}\Box}

\begin{document}

\title{Topological and algebraic structures on the ring of Fermat reals\date{}}

\author{Paolo Giordano\thanks{University of Vienna, Faculty of Mathematics, Nordbergstr. 15, A-1090, Austria. 
Supported by a L. Meitner FWF (Austria) grant (M1247-N13).}
\and Michael Kunzinger\thanks{University of Vienna, Faculty of Mathematics, Nordbergstr. 15, A-1090, Austria.
Supported by FWF research grants Y237 and P20525}}

\maketitle

\begin{center}
\emph{University of Vienna}
\par

\end{center}

\begin{abstract}
The ring of Fermat reals is an extension of the real field containing
nilpotent infinitesimals, and represents an alternative to Synthetic
Differential Geometry in classical logic. In the present paper, our
first aim is to study this ring from using standard topological and
algebraic structures. We present the Fermat topology, generated by a
complete pseudo-metric, and the omega topology, generated by a complete
metric. The first one is closely related to the differentiation of (non
standard) smooth functions defined on open sets of Fermat reals. The
second one is connected to the differentiation of smooth functions
defined on infinitesimal sets. Subsequently, we prove that every
(proper) ideal is a set of infinitesimals whose order is less than or
equal to some real number. Finally, we define and study roots of
infinitesimals. A computer implementation as well as an application to
infinitesimal Taylor formulas with fractional derivatives are presented.\\

\end{abstract}

\section{Introduction}
If mathematics is the language of nature, one can imagine that the
more results are discoverable and describable using a given part of
mathematics, the more faithfully that language will correspond to some
given part of nature. We can hence imagine that there is a sort of
weak isomorphism between that language and the corresponding part of nature.
Therefore, if two different languages are able to describe faithfully
the same part of nature, we can also think that these two languages
are, in some way, isomorphic to each other. So, because we are able to
use mathematical analysis and actual infinitesimals as languages to
describe nature, we can imagine that using the first one it would
be possible to obtain a rigorous and modern model of the informal
use of infinitesimals. If this idea is in some way correct, this should
actually be feasible without any non trivial background of mathematical
logic. The theory of Fermat reals represents a possible formalization
of this philosophical idea. Other possible approaches following this
line of thought are: Weil functors \cite{Kr-Mi2,Kr-Mi}, Levi-Civita
fields \cite{Lev1,Sha,Sh-Be}, Surreal numbers \cite{Con2,Ehrl},
geometries over a general base ring \cite{Ber}, or the ring of Colombeau
generalized numbers \cite{Col1,Col2,O}. Classical approaches requiring a non
trivial background of mathematical logic are Nonstandard Analysis
\cite{Rob66} and Synthetic Differential Geometry \cite{Koc,Lav,Mo-Re}.
In case the above philosophical idea sounds natural to the reader,
an open problem is to understand, from a mathematical, social or historical
point of view, why the latter theories, i.e. those requiring a non
trivial background of mathematical logic, seem more powerful than
the former ones.

The ring $\ER$ of Fermat reals can be defined and studied
using elementary calculus only (\cite{Gio10a}). It extends the field
$\R$ of real numbers and contains nilpotent infinitesimals, i.e.
$h\in\ER_{\ne0}$ such that $h^{n}=0$ for some $n\in\N_{>1}$.
The methodological thread followed in the development of the theory
of Fermat reals has always been guided by the necessity to obtain a good dialectic
between formal properties and their informal interpretations. Indeed,
the ring is totally ordered and geometrically representable (\cite{Gio10b,Gio10d,Gio09}),
to cite some examples.

Every Fermat real $x\in\ER$ can be written, in a unique way, as

\begin{equation}
x=\st{x}+\sum_{i=1}^{N}\alpha_{i}\cdot\diff{t}_{a_{i}},\label{eq:decomposition}
\end{equation}

\noindent where $\st{x}$, $\alpha_{i}$, $a_{i}\in\R$ are standard reals,
$a_{1}>a_{2}>\dots>a_{N}\ge1$, $\alpha_{i}\ne0$, and where $\diff{t}_{a}$
verifies the following properties

\begin{align}
\diff{t}_{a}\cdot\diff{t}_{b} & =\diff{t}_{\frac{ab}{a+b}}\nonumber \\
\left(\diff{t}_{a}\right)^{p} & =\diff{t}_{\frac{a}{p}}\quad\forall p\in\R_{\ge1}\label{eq:propertiesOfDiffInfinitesimals}\\
\diff{t}_{a} & =0\quad\forall a\in\R_{<1}.\nonumber 
\end{align}

The expression \eqref{eq:decomposition} is called the \emph{decomposition}
of $x$, and the real number $\st{x}$ its \emph{standard part}. The
number $a_{1}=:\omega(x)$ is called the \emph{order} of
$x$ and represents the greatest infinitesimal appearing in its decomposition.
In case $x\in\R,$ i.e. $x=\st{x}$, we set $\omega(x)=0$. We will
also use the notations $\omega_{i}(x):=a_{i}$ and $\st{x}_{i}:=\alpha_{i}$
for the $i$-th order and the $i$\emph{-th standard part} of $x$;
$\omega_{i}(x):=0$ if $x\in\R$. The order $\omega(-)$ has the following
natural properties

\begin{align*}
\omega(x+y) & =\max\left[\omega(x),\omega(y)\right]\\
\frac{1}{\omega(x\cdot y)} & =\frac{1}{\omega(x)}+\frac{1}{\omega(y)},
\end{align*}
whenever $x$, $y$ are infinitesimals such that $x+y\ne0$ or
$x\cdot y\ne0$, respectively.

Directly from \eqref{eq:decomposition} it is not hard to prove that
if $k\in\N_{>1}$, then $x^{k}=0$ iff $\omega(x)<k$. For $a\in\R_{\ge0}\cup\{\infty\}$,
the ideal
\[
D_{a}:=\left\{ x\in\ER\,|\,\st{x}=0,\ \omega(x)<a+1\right\} 
\]
plays a fundamental role in $k$-th order Taylor formulas with nilpotent
increments (so that the remainder is zero). Indeed, for $k\in\N_{\ge1}$
we have that $D_{k}=\left\{ x\in\ER\,|\, x^{k+1}=0\right\} $, and
any ordinary smooth function $f:A\freccia\R$ defined on an open set
$A$ of $\R^{d}$ can be extended to the set
\[
\ext{A}=\left\{ x\in\ER^{d}\,|\,\st{x}\in A\right\} ,
\]
obviously obtaining a true extension, i.e.\ the same values at $x\in A$.
The mentioned Taylor formula is therefore 
\[
\forall h\in D_{k}^{d}\pti f(x+h)=\sum_{\substack{j\in\N^{d}\\
|j|\le k}
}\frac{h^{j}}{j!}\cdot\frac{\partial^{|j|}f}{\partial x^{j}}(x),
\]
where $x\in A$ is a standard point, and $D_{k}^{d}=D_{k}\times\ptind^{d}\times D_{k}$.

It may seem difficult to work in a ring with zero divisors, but the
following properties permit to deal effectively with products of nilpotent
infinitesimals (typically appearing in multidimensional Taylor formulas)
and with cancellation laws:
\[
 h_1^{i_1}\cdot\ldots\cdot h_n^{i_n}=0\quad\iff\quad\sum_{k=1}^n\frac{i_k}{\omega(h_{k})}>1
\]

\[
x\text{ is invertible}\quad\iff\quad\st{x}\ne0
\]

\[
\left(\text{If }x\cdot r=x\cdot s\text{ in }\ER,\text{ where }r,s\in\R\text{ and }x\ne0\right)\then r=s.
\]
Finally, the ring $\ER$ is totally ordered, and the order relation
can be effectively decided, once again, starting from the decompositions:
let $x$, $y\in\ER$; if $\st{x}\ne\st{y}$, then

\begin{equation}
x<y\iff\st{x}<\st{y}.\label{eq:orderAndStdParts}
\end{equation}

\noindent Otherwise, if $\st{x}=\st{y}$, then

\begin{enumerate}
\item If $\omega(x)>\omega(y)$, then $x>y$ iff $\st{x_{1}}>0$, i.e. iff
	$x>0$ (from \eqref{eq:orderAndStdParts})
\item If $\omega(x)=\omega(y)$, then

	\begin{align*}
	\st{x_{1}}>\st{y_{1}} & \then x>y\\
	\st{x_{1}}<\st{y_{1}} & \then x<y.
	\end{align*}

\end{enumerate}

For example, $\diff{t}_{3}-3\diff{t}>\diff{t}>0$, and $\diff{t_{a}}<\diff{t_{b}}$
if $a<b$. This quick summary of some algebraic and order properties
of the ring $\ER$ of Fermat reals can be considered as a first step
toward its axiomatic description. A very simple model can already be
guessed from the properties \eqref{eq:decomposition} and \eqref{eq:propertiesOfDiffInfinitesimals}.
Indeed, we first introduce the ring $\R_{o}[t]$ of little-oh polynomials,
i.e. functions $x:\R_{\ge0}\freccia\R$ that can be written as $x(t)=r+\sum_{i=1}^{N}a_{i}\cdot t^{\alpha_{i}}+o(t)$,
as $t\to0^{+}$, where $N\in\N$, $r$, $a_{1},\dots,a_{N}\in\R$
and $\alpha_{1},\dots,\alpha_{N}\in\R_{\ge0}$. Then, in the ring
$\R_{o}[t]$, we define the equivalence relation $x\sim y$ iff $x(t)=y(t)+o(t)$,
for $t\to0^{+}$, and $\ER:=\R_{o}[t]/\sim$ is the related quotient
set. 

For a more complete presentation and study of the model, see \cite{Gio10a,Gio10d,Gio09}.

In \cite{Gio10e}, the Fermat-Reyes theorem, which is essential for
the development of differential calculus on $\ER$, is presented.
The Fermat-Reyes theorem states the existence and uniqueness of the
\emph{smooth incremental ratio} of every smooth function $f:\ext{U}\freccia\ER$,
that is existence and uniqueness of a function $r:\widetilde{\ext{U}}\freccia\ER$
satisfying
\[
f(x+h)=f(x)+h\cdot r(x,h)\quad\text{in}\quad\ER\quad\forall(x,h)\in\widetilde{\ext{U}},
\]
where $\widetilde{\ext{U}}:=\left\{ (x,h)\in\ER^{2}\,|\,\overrightarrow{[x,x+h]}\subseteq\ext{U}\right\} $
is called the thickening of $\ext{U}$ (\cite{Ab-Ma-Ra,Gio10e}).
Here, the function $f$ is more general than the extension from $U\subseteq\R$
to $\ext{U}\subseteq\ER$ of an ordinary smooth function defined on
$U$ and with values in $\R$. The function $f$ is a \emph{non standard
smooth} (or, more simply, \emph{smooth}) function, i.e., by definition,
$f$ can locally be written as 
\[
f(x)=\ext{\alpha}(p,x)\quad\forall x\in\ext{V}\cap\ext{U},
\]
where $\alpha\in\Cc^{\infty}(W\times V,\R^{\sf t})$ is an ordinary
smooth function defined on an open set of $\R^{\sf w}\times\R^{\sf v}$and
$p\in\ext{W}$ is a $\sf w$-dimensional parameter. The mentioned
topology is that generated by extension of open sets, i.e. by sets
of the form $\ext{U}$. For example, $f(x)=\diff{t}+x$ is (non standard)
smooth, but it is not the extension of an ordinary smooth function
because $f(0)=\diff{t}\notin\R$, whereas any extension takes $\R$
to itself. The Fermat-Reyes theorem is well framed in the cartesian
closed category $\ECInfty$ of Fermat spaces, as this permits
to develop a notion of smooth space and smooth function including
also infinite dimensional spaces, e.g. function spaces like $\ManInfty(M,N)$
or integral and differential operators. The definition of $\ECInfty$
is essentially a generalization of the notion of diffeological space
(\cite{Igl}), whose category $\CInfty$ is the domain of the Fermat
functor $\ext{(-)}:\CInfty\freccia\ECInfty$. This functor generalizes
the construction $\R\mapsto\ER$ and displays very good preservation
properties, closely related to intuitionistic logic. For more details
see \cite{Gio10e,Gio09}.

The main aim of the present paper is to study the ring of Fermat reals
using standard topological and algebraic structures. We will analyze interesting
metric structures deeply related to the development of smooth calculus
on $\ER$. We will characterize the ideals of $\ER$, confirming that
they are of a simple (non-pathological) nature, due to the initial
choice of the very well behaved family of little-oh polynomials. Finally
we will show that, in spite of the presence of infinitesimal numbers
$h\in\ER_{\ne0}$ such that $h^{2}=0$, we can define powers $h^{p}$
for every $p\in\R_{\ge0}$ and hence, we have arbitrary roots with
several good properties. This dialectic, between standard structures
and the new ring $\ER$, aims at presenting the theory of Fermat reals
to a general mathematical audience, hoping that this will contribute
to further studies of this interesting ring with infinitesimals.

\section{Metric structures}

The topology used to prove the above mentioned Fermat-Reyes theorem,
the key theorem for the development of differential calculus of smooth
functions defined on open sets, is that generated by extensions $\ext{U}$
of open sets $U\subseteq\R^{n}$. In this approach, a subset $A\subseteq\ER^{n}$
is open if it can be written as
\[
A=\bigcup\left\{ \ext{U}\subseteq A\,|\, U\text{ is open in }\R^{n}\right\} .
\]
We will call the resulting topology the \emph{Fermat topology}. In the present section,
we want to show that in the ring $\ER$ it is possible to define two
interesting (pseudo) metric structures, corresponding to two different
topologies. The first one is the Fermat topology, which can roughly be
described as the best topology for sets having a {}``sufficient amount
of standard points'', like, e.g., $\ext{U}$. This connection between
Fermat topology and standard points can be glimpsed by saying that the
monad of a standard real $r\in\R$, i.e.
\[
\mu(t):=\left\{ x\in\ER\,|\,\st{x}=r\right\} =\left\{ r+h\,|\, h\in D_{\infty}\right\} ,
\]
is the set of all the points which are limits of sequences with respect
to the Fermat topology (which is not Hausdorff).

However, in sets of infinitesimals, like the ideal $D_{a}$, there
is only one standard point, and indeed the best topology to study
this kind of sets is not the Fermat one. Therefore, we will define
a metric generating a finer topology, called the \emph{omega topology}.
When restricted to $\mu(r)$, the omega topology is naturally tied
with the equality $=_{k}$ up to $k$-th order infinitesimals, i.e.
$x=_{k}y$ iff $\st{x}=\st{y}$ and $\omega(x-y)\le k$ (see \cite{Gio10b,Gio09}
for some properties and applications of this notion). It is worth
noting that the equivalence relation $=_{k}$ is tied with differential
calculus of smooth functions defined on infinitesimal sets like $D_{a}$.
Indeed, in \cite{Gio09} it is proved that for a smooth function $f:D_{a}\freccia\ER$
there always exist $m_{1},\dots,m_{n}\in\ER$, $n:=[a]$ being the
integer part of $a$, such that
\[
f(h)=f(0)+\sum_{j=1}^{n}\frac{h^{j}}{j!}\cdot m_{j}\quad\forall h\in D_{a}.
\]
Moreover, these $m_{j}$ are unique up to $k_{j}$-th order infinitesimals
\[
\left(
	\forall h\in D_{a}:\ \sum_{j=1}^{n}\frac{h^{j}}{j!}\cdot m_{j}=\sum_{j=1}^{n}\frac{h^{j}}{j!}\cdot\bar{m}_{j}
\right)
\then m_{j}=_{k_{j}}\bar{m}_{j}\quad\forall j=1,\dots,n,
\]
where the $k_{j}$ are defined by
\[
\frac{1}{k_{j}}+\frac{j}{a+1}=1.
\]
This permits to define derivatives of smooth functions defined on
infinitesimal sets. Therefore, it is worth noting that two standard
metrics (i.e. with values in $\R$ and not in $\ER$) are strictly
related to the calculus of two different classes of smooth functions
on the ring of Fermat reals.

To motivate the definition of our metrics on $\ER$, we can say that:

\begin{itemize}
\item We want to measure the distance between $x$, $y\in\ER$ on the basis
of $\omega(x-y)$ and $|\st{x}-\st{y}|$ only.
\item We want to extend the classical metric on the reals $\R$.
\end{itemize}

\begin{defn}
Let $x$, $y\in\ER$, then

	\begin{enumerate}
	\item $d_{\text{F}}(x,y):=|\st{x}-\st{}y|$
	\item $d_{\omega}(x,y):=|\st{x}-\st{}y|+\omega(x-y).$
	\end{enumerate}

\end{defn}
Obviously, both $d_{\text{F}}$ and $d_{\omega}$ extend the usual
metric on $\R$; moreover, if $x$, $y\in\mu(r)$, then $d_{\omega}(x,y)\le k$
iff $x=_{k}y$. Finally, as we will see later, the idea to use also
the $i$-th orders $\omega_{i}(x-y)$ to define other metrics on $\ER$
is not a successful one.

\begin{thm}
The Fermat and the omega metrics verify the following properties:

\begin{enumerate}
\item \label{enu:dF_pseudometric}$d_{\text{F}}:\ER\times\ER\freccia\R$
is a pseudometric.
\item \label{enu:dOmega_metric}$d_{\omega}:\ER\times\ER\freccia\R$ is
a metric.
\item \label{enu:dOmega_finer_dF}The $d_{\omega}$-topology is finer than
the $d_{\text{F}}$-topology.
\item \label{enu:dF_generatesFermat}The topology generated by $d_{\text{F}}$
is the Fermat topology.
\item \label{enu:dF_dOmega_notEquivalent}$d_{\text{F}}$ and $d_{\omega}$
are not topologically equivalent.
\end{enumerate}
\end{thm}

\noindent\textbf{Proof:} The proof of \emph{\ref{enu:dF_pseudometric}} is
direct. Concerning \emph{\ref{enu:dOmega_metric}}, we have that $d_{\omega}(x,y)=0$
iff $\st{x}=\st{y}$ and $\omega(x-y)=0$. The order $\omega(x-y)$
is zero iff $x-y\in\R$, i.e. $x-y=:r\in\R$. From $\st{x}=\st{y}$
it follows that $r=0$ and hence the conclusion $x=y$. From the definition
of order $\omega(-)$, the property $\omega(-z)=\omega(z)$ follows,
and hence $d_{\omega}$ is symmetric. To prove
the triangle inequality, we introduce the following lemma, which is
a generalization of an analogous result already proved for $x$, $y\in D_{\infty}$
only (see e.g. \cite{Gio10a}).

\begin{lem}
\label{lem:orderAndSum}Let $x$, $y\in\ER$, then

\begin{enumerate}
\item \label{enu:LmOrderSum_notReal}$x+y\notin\R\then\omega(x+y)=\max\left[\omega(x),\omega(y)\right]$
\item \label{enu:enu:LmOrderSum_Real}$x+y\in\R\then\omega(x+y)=0$.
\end{enumerate}

Therefore
\[
\omega(x-y)\le\omega(x-z)+\omega(z-y)\quad\forall z\in\ER.
\]

\end{lem}
\noindent Using this lemma we have

\begin{align*}
d_{\omega}(x,y) & =|\st{x}-\st{y}|+\omega(x-y)\le\\
 & \le|\st{x}-\st{z}|+|\st{z}-\st{y}|+\omega(x-z)+\omega(z-y)=\\
 & =d_{\omega}(x,z)+d_{\omega}(z,y).
\end{align*}

Property \emph{\ref{enu:dOmega_finer_dF}} follows directly from the
inequality $d_{\text{F}}(x,y)\le d_{\omega}(x,y)$.

\noindent To prove \emph{\ref{enu:dF_generatesFermat}} let us firstly
consider a $d_{\text{F}}$-open set $\mathcal{U}$. Then, for every
$x\in\mathcal{U}$ we can find $r\in\R_{>0}$ such that

\begin{equation}
B_{r}(x;d_{\text{F}})=\left\{ y\in\ER\,|\,|\st{x}-\st{y}|<r\right\} \subseteq\mathcal{U}.\label{eq:ball_mathcal_U}
\end{equation}

Let $U:=(\st{x}-r,\st{x}+r)_{\R}=\{s\in\R\,|\,|\st{x}-s|<r\}$, then
for every $y\in\ext{U}$ we have $\st{y}\in U$ and hence $y\in\mathcal{U}$
from \eqref{eq:ball_mathcal_U}. Therefore, $x\in\ext{U}\subseteq\mathcal{U}$,
that is $\mathcal{U}$ is also open in the Fermat topology. Vice versa,
let $\mathcal{U}$ be an open set in the Fermat topology, then for
every $x\in\mathcal{U}$ we can find an open set $U$ of $\R$ such
that $x\in\ext{U}\subseteq\mathcal{U}$. Therefore, $\st{x}\in U$
and

\begin{equation}
(\st{x}-r,\st{x}+r)_{\R}\subseteq U\label{eq:stdInterval_U}
\end{equation}
for some $r\in\R_{>0}$. So, for every $y\in B_{r}(x;d_{\text{F}})$
we have $|\st{x}-\st{y}|<r$ and $\st{y}\in U$ from \eqref{eq:stdInterval_U}.
This implies $y\in\ext{U}\subseteq\mathcal{U}$ and proves that $B_{r}(x;d_{\text{F}})\subseteq\mathcal{U}$.

To prove \emph{\ref{enu:dF_dOmega_notEquivalent}} we can consider
$\mathcal{U}=B_{r}(0;d_{\omega})$, where $r>1$. We want to show
that every ball $B_{s}(0;d_{\text{F}})$ is not contained in $\mathcal{U}$,
that is
\[
\exists x\in\ER:\ |\st{x}|<s\ ,\ |\st{x}|+\omega(x)\ge r.
\]
To prove this it suffices to show that $\omega(x)\ge r$ for some
infinitesimal $x\in D_{\infty}$, which is trivially true: we can
take, e.g., $x=\diff{t}_{r}$ whose order is $\omega(\diff{t}_{r})=r$
because $r>1$.$\qed$ 

\noindent \textbf{Proof of Lemma \ref{lem:orderAndSum}:} Let $\delta x:=x-\st{x}$,
$\delta y:=y-\st{y}$ be the infinitesimal parts of $x$, $y$, so
that, directly from the definition of order, we have $\omega(x+y)=\omega(\delta x+\delta y)$.
If $\delta x+\delta y=0$, then $\omega(x+y)=\omega(0)=0$. Otherwise,
$\delta x+\delta y\ne0$ and from Theorem 12 of \cite{Gio10a} we
have $\omega(\delta x+\delta y)=\max\left[\omega(\delta x),\omega(\delta y)\right]=\max\left[\omega(x),\omega(y)\right].$
Finally, it suffices to note that
\[
\delta x+\delta y=0\iff x+y\in\R
\]
to prove both properties \emph{\ref{enu:LmOrderSum_notReal}} and
\emph{\ref{enu:enu:LmOrderSum_Real}}.

To prove the stated inequality, we note that if $x+y\in\R$,
then $\omega(x+y)=0\le\omega(x)+\omega(y)$ because $\omega(-)\ge0$.
If $x+y\notin\R$, then either $\omega(x+y)=\omega(x)\le\omega(x)+\omega(y)$
or $\omega(x+y)=\omega(y)\le\omega(x)+\omega(y)$. From this, the
conclusion follows, because $\omega(x+y)=\omega\left[(x-z)+(z-y)\right]$.$\qed$

\begin{defn}
We will call $\omega$\emph{-topology} the topology generated by the
metric $d_{\omega}$. It can also be called the topology of the \emph{order
function} (to distinguish it from the topology of the order relation).
\end{defn}

\subsubsection*{Is it possible to generalize using higher orders $\omega_{i}$?}

It is very natural to try a generalization of the metric $d_{\omega}$
considering, e.g., also the information given by $\omega_{2}(x-y)$:
\[
d_{2}(x,y):=|\st{x}-\st{y}|+\omega(x-y)+\omega_{2}(x-y).
\]
However, an immediate problem is that higher orders $\omega_{i}(-)$,
$i>1$, do not share the good properties of $\omega(-)$. For example:

\begin{itemize}
\item If $x=\diff{t}$ and $y=\diff{t}_{3}+\diff{t}_{2}$, then $\omega_{2}(x+y)=\omega_{2}(y)$,
but if $y=\diff{t}_{2}$, then $\omega_{2}(x+y)=\omega(x)$.
\item If $\omega(x)=\omega(y)$ and $\st{x}_{1}+\st{y}_{1}\ne0$, then $\omega_{2}(x+y)=\max\left[\omega_{2}(x),\omega_{2}(y)\right]$,
but if $\st{x}_{1}+\st{y}_{1}=0$, $\omega_{2}(x)=\omega_{2}(y)$
and $\st{x}_{2}+\st{y}_{2}\ne0$, then $\omega_{2}(x+y)=\max\left[\omega_{3}(x),\omega_{3}(y)\right]$.
\item If $x=\diff{t}_{6}+2\diff{t}$, $y=\diff{t}_{7}+\diff{t}$ and $z=3\diff{t}$,
then $\omega_{2}(x-z)+\omega_{2}(z-y)<\omega_{2}(x-y)$. Therefore,
the inequality of Lemma \ref{lem:orderAndSum} cannot be proved for
higher orders.
\end{itemize}

We will solve this problem with the following result.

\begin{prop}
Let $i\in\N_{>1}$, and suppose that
\[
d_{i}(x,y):=|\st{x}-\st{y}|+\sum_{j=1}^{i}\omega_{j}(x-y)
\]
verifies the triangle inequality. Then $d_{i}$ and $d_{\omega}$
are equivalent.
\end{prop}

\begin{rem}
In this statement, we mean $\omega_{j}(r):=0$ for $r\in\R$ and $\omega_{j}(x):=0$
if $j>N$, where $N$ is the number of summands in the decomposition
of $x$.
\end{rem}

\noindent \textbf{Proof:} It is readily verified that the $d_{i}$-topology
is finer than the $\omega$-topology. To prove the converse, let us
first consider $y\in B_{r}(x;d_{\omega})$, with $r\le1$. Then
$d_{\omega}(x,y)=|\st{x}-\st{y}|+\omega(x-y)<r\le1$. Therefore, $\omega(x-y)<1$
and hence $\omega(x-y)=0$. This means $x-y\in\R$ and hence $\omega_{j}(x-y)=0$
and $d_{i}(x,y)=d_{\omega}(x,y)<r$. Therefore, $B_{r}(x;d_{\omega})\subseteq B_{r}(x;d_{i})$
if $r\le1$. Finally, for a generic ball $B_{s}(x;d_{i})$, take $n\in\N_{>0}$
such that $\frac{s}{n}\le1$, then $B_{s}(x;d_{i})\supseteq B_{s/n}(x;d_{i})\supseteq B_{s/n}(x;d_{\omega})$.$\qed$

\subsection{Fermat and $\omega$-completeness of $\ER$}

The proof of the following result follows directly from the equality
$d_{\text{F}}(s_{n},s_{m})=|\st{s_{n}}-\st{s_{m}}|$.

\begin{thm}
With respect to the Fermat metric $d_F$, the ring $\ER$ is complete. In particular,
if $(s_{n})_{n}$ is a Cauchy sequence with respect to $d_{\text{F}}$,
then $\left(\st{s_{n}}\right)_{n}$ is a standard Cauchy sequence
of $\R$. Let $r\in\R$ be its limit, then
\[
\forall x\in\mu(x):\ (s_{n})_{n}\text{ converges to }x\text{ with respect to }d_{\text{F}}.
\]
\end{thm}

\noindent Before studying the $\omega$-topology, we want to understand
better the intuition underlying this metric, because it is strictly
related to nilpotency of every infinitesimal of $\ER$. Let us start
with an example:
\[
s_{n}:=\frac{1}{n}+\diff{t}_{\frac{N}{n}}\quad\forall n\in\N_{>0},
\]
where $N\in\N_{>0}$. We have $d_{\omega}(s_{n},0)=\frac{1}{n}+\omega\left(\frac{1}{N}+\diff{t}_{\frac{N}{n}}-0\right)=\frac{1}{n}+\frac{N}{n}\to0$
as $n\to+\infty$. Therefore, $(s_{n})_{n}$ converges to 0 in the
$\omega$-topology. However, exactly because of nilpotency, we also
have
\[
\forall n>N:\ \diff{t}_{\frac{N}{n}}=0\quad\text{hence}\quad s_{n}=\frac{1}{n}\in\R.
\]
More generally, if $d_{\omega}(s_{n},x)\to0$, then $|\st{s_{n}}-\st{x}|\to0$,
but also $\omega(s_{n}-x)\to0$. This means that the order $\omega(s_{n}-x)$
goes through smaller and smaller infinitesimals. However, because
of nilpotency, infinitesimals of $\ER$ cannot have order less than
1, and hence the order $\omega(s_{n}-x)$ must collapse from 1 to
0. The following theorems will clarify this intuition.

\begin{thm}
Let $x\in\ER$ and $s\in(0,1]_{\R}$, then
\[
B_{s}(x;d_{\omega})=(-s,s)_{\R}+\{x\},
\]
that is, every $y\in B_{s}(x;d_{\omega})$ can be written as $y=r+x$,
with $r\in\R$, $-s<r<s$.
\end{thm}

\noindent \textbf{Proof:} For $y\in B_{s}(x;d_{\omega})$ we have
$|\st{y}-\st{x}|<s$ and $\omega(y-x)<s\le1$, therefore $\omega(y-x)=0$
and $y-x=:r\in\R$. Moreover, $r=\st{y}-\st{x}$ so that $|r|<s$,
that is $r\in(-s,s)_{\R}$. Vice versa, if $y\in(-s,s)_{\R}+\{x\}$,
then $y=r+x$, with $|r|<s$. So $|\st{y}-\st{x}|=|r|<s$ and $\omega(y-x)=\omega(r)=0$
and $d_{\omega}(y,x)<s$.$\qed$

The following theorem characterizes $\omega$-convergent sequences
formalizing the intuition presented above.

\begin{thm}
Let $(s_{n})_{n}$ be a sequence of $\ER$, then we have that

\begin{equation}
\lim_{n\to+\infty}s_{n}=x\in\ER\label{eq:limitCharactOmegaConvergence}
\end{equation}
with respect to the omega topology if and only if the following conditions hold

\begin{enumerate}
\item $\lim_{n\to+\infty}\st{s_{n}}=\st{x}$ in $\R$
\item \label{enu:definitelyConstant}The sequence of infinitesimal parts
is eventually constant and equal to $\delta x$, i.e.
\[
\exists N\in\N:\ \forall n\ge N:\ \delta s_{n}=\delta x.
\]
(where we recall that $\delta y:=y-\st{y}$.)
\end{enumerate}
\end{thm}

\noindent\textbf{Proof:} We only have to prove \emph{\ref{enu:definitelyConstant}},
as the rest of the proof is immediate. Because of \eqref{eq:limitCharactOmegaConvergence},
we have that
\[
\lim_{n\to+\infty}\omega(s_{n}-x)=0=\lim_{n\to+\infty}\omega\left(\delta s_{n}-\delta x\right).
\]
Therefore $\omega\left(\delta s_{n}-\delta x\right)<1$ for $n\ge N$,
and hence $\omega\left(\delta s_{n}-\delta x\right)=0$. This means
that $\delta s_{n}-\delta x\in\R$ and hence $\delta s_{n}-\delta x=0$
because $\delta s_{n}$, $\delta x\in D_{\infty}$.$\qed$

\begin{example}
$s_{n}=\frac{1}{n}+\diff{t}_{\frac{N}{n}}+\diff{t}_{2}\to\diff{t}_{2}$,
whereas $s_{n}=\diff{t}_{1+\frac{1}{n}}$ is not an $\omega$-convergent
sequence.
\end{example}

\begin{thm}
The ring $\ER$ is complete with respect to the metric $d_\omega$.
\end{thm}

\noindent \textbf{Proof:} Let $(s_{n})_{n}$ be an $\omega$-Cauchy
sequence, then
\[
\lim_{\substack{n\to+\infty\\
m\to+\infty}
}|\st{s_{n}}-\st{s_{m}}|=0,
\]
so that the $\lim_{n\to+\infty}\st{s_{n}}=:r\in\R$ exists. Moreover,
we have
\[
\lim_{\substack{n\to+\infty\\
m\to+\infty}
}\omega(s_{m}-s_{n})=0=\lim_{\substack{n\to+\infty\\
m\to+\infty}
}\omega(\delta s_{m}-\delta s_{n}).
\]
Like in the previous proof, we have $\delta s_{n}=\delta s_{m}=\delta s_{N}$
for $n$, $m\ge N$. Therefore, setting $x:=r+\delta s_{N}\in\ER$
we have $\lim_{n\to+\infty}s_{n}=x$.$\qed$

\subsubsection{The omega metric on $D_{\infty}$ can be defined by a pseudovaluation}

In this section we want to prove that the restriction of $d_{\omega}$
to the set of all the infinitesimals $D_{\infty}$
\[
d_{\omega}(h,k)=\omega(h-k)\quad\forall h,k\in D_{\infty}
\]
is induced by a pseudovaluation. We will see that
the same idea doesn't work outside $D_{\infty}$. These notes have
also the aim to fix some small error made on the same topic in \cite{Gio10a}.

For completeness, we start from the definition of pseu\-do\-va\-lua\-tion
on a generic ring with values in $\R\cup\{+\infty\}$.

\begin{defn}
\label{def:pseudovaluation}Let $A$ be a ring, then we say that $v:A\freccia\R\cup\{+\infty\}$
is a pseudovaluation if the following conditions hold:

\begin{enumerate}
\item $v:A\freccia\R\cup\{+\infty\}$
\item $\forall x\in A:\ v(x)=+\infty\iff x=0$
\item \label{enu:productAndPseudoVal}$v(x\cdot y)\ge v(x)+v(y)$
\item \label{enu:minPseudval}$v(x+y)\ge\min\left[v(x),v(y)\right]$
\item $v(x-y)=v(y-x)$
\end{enumerate}

Obviously, the last three conditions are supposed to hold for every $x$,
$y\in A$.
\end{defn}

\begin{rem}\ 
\begin{enumerate}
\item We assume the usual rules about the relationships between $+\infty$
and order or sum on $\R$: $+\infty+r=+\infty>r$ for every $r\in\R$.
\item The essential difference between a pseudovaluation and a valuation
is property \emph{\ref{enu:productAndPseudoVal}}, where an inequality
replaces an equality. Indeed, it is not hard to prove that equality
cannot hold in any ring with zero divisors, like $\ER$.
\end{enumerate}
\end{rem}

To motivate the necessity of our definition of pseudovaluation, we
anticipate the following lemma, with which we can treat terms of the
form $\omega(x\cdot y)$.

\begin{lem}
Let $x$, $y\in\ER$, then
\[
\omega(x\cdot y)=
\begin{cases}
\omega(x) & \text{ if }x\in D_{\infty},\ y\notin D_{\infty}\\
\left[\frac{1}{\omega(x)}+\frac{1}{\omega(y)}\right]^{-1} & \text{ if }x,y\in D_{\infty}\setminus\{0\}
\end{cases}
\]
Moreover, if $x$, $y\notin\R$, then 
\[
\omega(x\cdot y)\le\omega(x)\cdot\omega(y).
\]
\end{lem}

\noindent \textbf{Proof:} Using the infinitesimal parts, we can write
$x\cdot y=\st{x}\cdot\st{y}+\st{x}\cdot\delta y+\st{y}\cdot\delta x+\delta x\cdot\delta y$,
so that, by the definition of order, the equality $\omega(x\cdot y)=\omega(\st{x}\cdot\delta y+\st{y}\cdot\delta x+\delta x\cdot\delta y)$
follows.

In the last case of the statement, i.e. $x$, $y\in D_{\infty}\setminus\{0\}$,
we have

\begin{align*}
\omega(x\cdot y) & =\omega(\delta x\cdot\delta y)=\left[\frac{1}{\omega(x)}+\frac{1}{\omega(y)}\right]^{-1}=\frac{\omega(x)\cdot\omega(y)}{\omega(x)+\omega(y)}\le\\
 & \le\frac{1}{2}\omega(x)\cdot\omega(y)\le\omega(x)\cdot\omega(y)
\end{align*}

Moreover, let us also note that $\frac{1}{\omega(x\cdot y)}=\frac{1}{\omega(x)}+\frac{1}{\omega(y)}\ge\frac{1}{\omega(x)}$,
and hence

\begin{equation}
\omega(x\cdot y)\le\omega(x).\label{eq:inequalityOfOrdersForTwoInfinitesimals}
\end{equation}

Below we will use this inequality for suitable infinitesimals $x$
and $y$.

If $x\in D_{\infty}$ and $y\notin D_{\infty}$, the cases $\delta x=0$
or $\delta x\cdot\delta y=0$ are immediate. Otherwise, we have
\begin{align}
\omega(x\cdot y) & =\omega(\st{x}\cdot\delta y+\st{y}\cdot\delta x+\delta x\cdot\delta y)=\omega(\st{y}\cdot\delta x+\delta x\cdot\delta y)=\nonumber \\
 & =\max\left[\omega(\st{y}\cdot\delta x),\omega(\delta x\cdot\delta y)\right].\label{eq:2_secondCase}
\end{align}
Now, we can apply \eqref{eq:inequalityOfOrdersForTwoInfinitesimals}
to the product $\delta x\cdot\delta y$, obtaining

\begin{equation}
\omega(\delta x\cdot\delta y)\le\omega(\delta x)=\omega(\st{y}\cdot\delta x).\label{eq:3_secondCase}
\end{equation}

Note that the last equality is due to the fact that $y\notin D_{\infty}$,
so that $\st{y}\ne0$. Therefore, from \eqref{eq:2_secondCase} and
\eqref{eq:3_secondCase} we obtain $\omega(x\cdot y)=\omega(\st{y}\cdot\delta x)=\omega(\delta x)=\omega(x)$.
Finally, $\omega(x)\le\omega(x)\cdot\omega(y)$ if $y\notin\R$ because
in that case $\omega(y)\ge1$. 
$\qed$

\begin{rem}
In the statement of the previous theorem, the case where $x$, $y\notin D_{\infty}$
is not included. This is done to avoid an overcomplicated
statement. Indeed, we have several sub-cases:

\begin{enumerate}
\item If $\st{x}\cdot\delta y+\st{y}\cdot\delta x+\delta x\cdot\delta y=0\then\omega(x\cdot y)=0$
\item If $\st{x}\cdot\delta y+\st{y}\cdot\delta x=0$ and $\delta x\cdot\delta y\ne0\then\omega(x\cdot y)=\left[\frac{1}{\omega(x)}+\frac{1}{\omega(x)}\right]^{-1}$
\item If $\st{x}\cdot\delta y+\st{y}\cdot\delta x\ne0\then\omega(x\cdot y)=\max\left[\omega(x),\omega(y)\right]$.
\end{enumerate}

\end{rem}
The idea to define a pseudovaluation is to derive the property $v(x\cdot y)\ge v(x)+v(y)$
from $\omega(x\cdot y)\le\omega(x)\cdot\omega(y)$, so that the natural
try is the following

\begin{defn}
For every infinitesimal $h\in D_{\infty}$, we define
\[
v(h):=-\log\left[\omega(h)\right],
\]
where we use the convention that $-\log(0):=+\infty$.
\end{defn}

\begin{rem}\ 

\begin{enumerate}
\item The metric associated to $v$ is $e^{-v(h-k)}=\omega(h-k)=d_{\omega}(h,k)$.
\item If we define $v(x):=-\log\left[|\st{x}|+\omega(x)\right]$ for every
$x\in\ER$, we do not obtain a pseudovaluation. Indeed, we have e.g.
$v(3+5)<\min\left[v(3),v(5)\right]$.
\end{enumerate}

\end{rem}

\begin{thm}
$v:D_{\infty}\freccia\R\cup\{+\infty\}$ is a pseudovaluation on the
subring (ideal) $D_{\infty}$ of all the infinitesimals.
\end{thm}
\noindent\textbf{Proof:} We only have to prove property \ref{enu:minPseudval}
of Definition \ref{def:pseudovaluation}, as the others are immediate.
We can suppose $x+y\ne0$, because otherwise the proof is obvious.
Since $x$, $y\in D_{\infty}$ and $x+y\ne0$, we have that $\omega(x+y)=\max\left[\omega(x),\omega(y)\right]$.
We will proceed in the case $\omega(x)\ge\omega(y)$, the opposite
being analogous. Therefore, $\omega(x+y)=\omega(x)$ and $v(x+y)=-\log\left[\omega(x+y)\right]=v(x)$.
>From $\omega(x)\ge\omega(y)$ we have $-\log\omega(x)\le-\log\omega(y)$,
so that $v(x)\le v(y)$ and $\min\left[v(x),v(y)\right]=v(x)=v(x+y)$,
which is our conclusion.$\qed$

\section{Ideals and their characterization}

In this section, we want to study the ring of Fermat reals from the
point of view of some standard algebraic structure. 
To begin with, we note a few elementary algebraic properties of $\ER$:
\begin{itemize}
\item There are no nontrivial idempotents in $\ER$.
\item $\ER$ is not reduced.
\item $x\in \ER$ is a zero divisor iff $x\in D_\infty$ iff $x$ is non-invertible.
\item $\ER$ is an exchange ring (i.e., for each $x\in \ER$ there exists an
idempotent such that $x+e$ is invertible).
\item $\ER$ is an $l$-ring (lattice ordered ring), as well as a normal $f$-ring
(cf.\ \cite{BKW77}).
\end{itemize}

Next we want to study the ideals of $\ER$. We will see that,
as opposed to other rings containing infinitesimals
(see, e.g., \cite{Ver10}) the ideals in $\ER$ can be exhaustively described. Of
course, this is essentially due to the very simple family of little-oh
polynomials used as representatives of new numbers in $\ER$. We start by
proving that the only maximal ideal is the set $D_{\infty}$ of all
the infinitesimals.
The idea to consider only ``well behaved'' functions $x:\R_{\ge0}\freccia\R$
(the little-oh polynomials), in the definition of the ring $\ER$,
is tied with the fact that Fermat reals do not represent
a new foundation for the entire calculus.
Indeed, our aim is only to extend ordinary smooth functions, 
so that it suffices to evaluate them on `well-behaved numbers'.

The situation is entirely different in NSA, which
aims to be a new, independent foundation of the whole calculus. For
example, suppose we want to prove that ordinary continuity of a function $f:\R\freccia\R$
at $x_{0}\in\R$ is equivalent to
\begin{equation}
\forall x\in\hyperR:\ x\simeq x_{0}\then{^{*}f}(x)\simeq{^{*}f}(x_{0}),\label{eq:NSA_continuitiy}
\end{equation}
where $x\simeq y$ means that $x-y$ is infinitesimal. However, \eqref{eq:NSA_continuitiy}
is nothing more than the continuity of the function $f$ stated using
sequences, i.e.

\begin{equation}
\forall x\in\R^{\N}:\ \exists\lim_{n\to+\infty}x_{n}=x_{0}\then\exists\lim_{n\to+\infty}f(x_{n})=f(x_{0}).\label{eq:continuityBySequences}
\end{equation}

This is equivalent to ordinary continuity only if it is stated for
every sequence $x\in\R^{\N}$, as is obvious from the corresponding
proof.

\begin{lem}
Let $I$ be a proper ideal of the ring $\ER$, then $I\subseteq D_{\infty}$.
\end{lem}

\noindent\textbf{Proof:} Let $x\in I$ and suppose, by contradiction, that
$\st{x}\ne0$, then $x$ would be invertible and for every $y\in\ER$
we could write $y=y\cdot x^{-1}\cdot x$. By hypothesis $x\in I$, which
is an ideal, so we would have $y\in I$, that is $I=\ER$ which is
impossible because $I$ is proper by hypothesis.$\qed$

Directly from the decomposition of every Fermat real it follows
that
\[
\ER/D_{\infty}\simeq\R
\]
and hence $D_{\infty}$ is a maximal ideal. With the following result,
we prove that in fact it is the only one.

\begin{thm}
Let $I$ be a proper maximal ideal of the ring $\ER$, then $I=D_{\infty}$.
\end{thm}

\noindent \textbf{Proof:} By the previous lemma, we have that $I\subseteq D_{\infty}$.
The ring $\ER$ is commutative and with unity, so $\ER/I$ is a field
because, by hypothesis, $I$ is maximal. Take $x\in D_{\infty}$.
We have two cases: either $x+I=I$ or $x+I\ne I$. In the first one, we have
$x+0=x\in I$. In the second one, $x+I$ would be invertible in the
field $\ER/I$, so, for some $y\in\ER$, we can write
\[
\left(y+I\right)\cdot\left(x+I\right)=1+I\quad\text{i.e.}\quad yx+I=1+I,
\]
that is $yx+i=1+j$ for some $i$, $j\in I$. Taking the standard
parts in this equality we obtain $\st{y}\cdot\st{x}+\st{i}=1+\st{j}$.
>From $I\subseteq D_{\infty}$ we deduce that $\st{i}=\st{j}=0$ and,
therefore, that $\st{y}\cdot\st{x}=1$, that is $\st{x}\ne0$, which
is impossible because $x\in D_{\infty}$.$\qed$

This proof can be easily generalized to the following

\begin{thm}
Let $R$, $F$ be two commutative rings with unity, with $F$ non
trivial. Moreover, let $s:R\freccia F$ be a ring morphism and set
$D:=\ker(s)$. Finally, let us suppose that
\[
\forall x\in R:\ x\notin D\then x\text{ is invertible in }R.
\]
Then $D$ is the only maximal ideal of $R$.
\end{thm}

In our case $s=\st{(-)}$ is the standard part map.

Finally, we prove that every proper ideal is either of the form
\[
D_{a}=\left\{ x\in D_{\infty}\,|\,\omega(x)<a+1\right\} \quad\text{for some}\quad a\in\R_{\ge1}\cup\{\infty\},
\]
or of the form
\[
I_{a}:=\left\{ x\in D_{\infty}\,|\,\omega(x)\le a\right\} \quad\text{for some}\quad a\in\R_{\ge1}.
\]
As we mentioned above, the first type of ideal is used in infinitesimal
Taylor formulas, whereas the second is used in the study of the
equivalence relation $=_{a}$ of equality up to $a$-th order infinitesimals.

This characterization is tied to the possibility to solve in $\ER$
the following class of linear equations.

\begin{thm}
\label{thm:linearEquations}If $a$, $b$, $c\in\ER$ and $a<c<a+b$,
then
\[
\exists x\in\ER:\ a+x\cdot b=c
\]
\end{thm}

\noindent For the proof of this theorem, see \cite{Gio10e,Gio09}.
Let us note that we cannot have uniqueness of solutions, due to nilpotency.
For example, if $a=0$, $c=\diff{t}_{2}+\diff{t}$ and $b=\diff{t}_{3}$,
then $x=\diff{t}_{6}+\diff{t}_{3/2}$ is a solution of $a + x\cdot b=c$,
but $x+\diff{t}$ is another solution. Moreover, let us note that
this theorem is not in contradiction with the non Archimedean property
of $\ER$ (let $a=0$ and $b\in D_{\infty})$ because of the inequalities
that $c$ must verify for a solution to exist.

Using this result, we can prove the desired characterization:

\begin{thm}
Let $J$ be a proper ideal of $\ER$, and $b\in J$. Moreover, set
\[
O(J):=\left\{ \omega(j)\in\R_{\ge0}\,|\, j\in J\right\} \quad,\quad a:=\sup O(J).
\]
Then

\begin{enumerate}
\item \label{enu:1_charIdeals}$\forall c\in\ER:\ -b<c<b\then c\in J$
\item \label{enu:2_charIdeals}$a\in O(J)\then J=I_{a}$
\item \label{enu:3_charIdeals}$a\notin O(J)\then a\ge1\quad\text{and}\quad J=D_{a-1}.$
\end{enumerate}

\end{thm}

\noindent \textbf{Proof:} To prove \emph{\ref{enu:1_charIdeals}},
set $a=0$ in Theorem \ref{thm:linearEquations}. Then since
$-b<c<b$, we can distinguish two cases (we recall that the order
relation in $\ER$ is total). If $c\ge0$, then $0\le c<b$ and we
can hence solve the equation $x\cdot b=c$ and, therefore, $c\in J$
because $b\in J$. Otherwise, $c<0$ and so $0<-c<b$. We solve the
equation $x\cdot b=-c$, that is $(-x)\cdot b=c$ so that $c\in J$
again.

To prove \emph{\ref{enu:2_charIdeals}} we first note that
$J\subseteq I_{a}$ by the definition of $a$. Vice versa, let $i\in I_{a}$,
i.e. $\omega(i)\le a$. By hypothesis, $a\in O(J)$, so that we can
write $a=\omega(j)$ for some $j\in J$. We can suppose $j>0$ because,
otherwise, $0<-j\in J$ and $a=\omega(j)=\omega(-j)$. We distinguish
two cases. If $\omega(i)<a=\omega(j)$, then $i<j$ by the properties
of the order relation we mentioned in the introduction (see Theorem
4.2.6 in \cite{Gio09}). On the other hand, we also have that $\omega(-i)=\omega(i)<\omega(j)$
and hence $-i<j$ because $j>0$. Therefore, $-j<i<j$, and the conclusion
$i\in J$ follows from \emph{\ref{enu:1_charIdeals}}. Let us note
that, in general, we have just proved that

\begin{equation}
\forall i,j\in\ER:\ j>0\ ,\ \omega(i)<\omega(j)\then-j<i<j.\label{eq:0_charIdeals}
\end{equation}
In the second case, we suppose that $\omega(i)=a$ so that, by the
decompositions of $i$, $j$ and for suitable $\alpha$, $\beta\in\R_{\ne0}$
and $h$, $k\in D_{\infty}$, we can write

\begin{align}
j & =\alpha\cdot\diff{t}_{a}+h\quad,\quad\omega(h)<a\label{eq:1_charIdeals}\\
i & =\beta\cdot\diff{t}_{a}+k\quad,\quad\omega(k)<a.\label{eq:2_charIdeals}
\end{align}
Therefore, from \eqref{eq:1_charIdeals} and \eqref{eq:0_charIdeals}
it follows that $-j<h<j$ and hence $h\in J$ from property \emph{\ref{enu:1_charIdeals}}.
So $\alpha\cdot\diff{t}_{a}=j-h\in J$ and $\diff{t}_{a}\in J$ because
$\alpha\ne0$. Hence, we also have that $\beta\cdot\diff{t}_{a}\in J$.
Finally, $\omega(k)<a=\omega(j)$ so that $-j<k<j$ from \eqref{eq:0_charIdeals}
and $k\in J$ from property \ref{enu:1_charIdeals}. We have proved
that $\beta\cdot\diff{t}_{a}$, $k\in J$, so $i=\beta\cdot\diff{t}_{a}+k\in J$,
which is our conclusion.

Finally, to prove \emph{\ref{enu:3_charIdeals}} we first note
that
\[
D_{a-1}=\left\{ x\in D_{\infty}\,|\,\omega(x)<a\right\} ,
\]
where we use the conventions $\infty\pm1=\infty$. If $a\notin O(J)$
then we have $\omega(j)<a$ for every $j\in J$, and therefore $J\subseteq D_{a-1}$,
considering also that $J\subseteq D_{\infty}$.

Vice versa, if $i\in D_{a-1}$, then we can find $j\in J$ such that
$\omega(i)<\omega(j)<a$ because $a=\sup O(J)$. As usual, we can
suppose $j>0$. From \eqref{eq:0_charIdeals}, it follows $-j<i<j$
and hence $i\in J$ by property \emph{\ref{enu:1_charIdeals}}. To
finish, let us note that because $\omega(j)\ge1$ or $\omega(j)=0$
we necessarily have that $a=0$ or $a\ge1$. However, the first possibility
would imply $J=\{0\}$ and hence $a\in O(J)$, which is impossible
by hypothesis.$\qed$ 

\section{Roots of infinitesimals}

In the ring of Fermat reals $\ER$, the existence of non zero nilsquare
elements:

\begin{equation}
h\ne0\quad,\quad h^{2}=0,\label{eq:nonZeroNilsquare}
\end{equation}
is incompatible with the existence of a square root and of an absolute
value with the usual properties. In other words, if we want to define
roots of infinitesimals, we have to avoid from \eqref{eq:nonZeroNilsquare}
the following inference:
\begin{align*}
h^{2} & =0\quad\text{therefore}\quad\sqrt{h^{2}}=\sqrt{0}=0\\
\sqrt{h^{2}} & =|h|=0\quad\text{therefore}\quad h=0.
\end{align*}

We recall that only smooth functions $f:\R\freccia\R$ can be extended
to $\ER$. In particular:

\begin{itemize}
\item Because they are locally Lipschitz, these functions verify
\[
\forall x,y\in\ER:\ x=y\text{ in }\ER\then f\circ x=f\circ y\text{ in }\ER.
\]
\item Because they are smooth, they take little-oh polynomials into little-oh
polynomials:
\[
\forall x\in\R_{o}[t]:\ f\circ x\in\R_{o}[t].
\]
\end{itemize}

\noindent It is hence natural to expect some problems extending, e.g.,
the square root to the whole of $\ER$.

The first natural solution is to extend the roots only where they
are smooth, i.e. on $\R_{\ne0}$. This is equivalent to defining the
roots only for invertible Fermat reals (and positive in case of even
roots or irrational powers). For details about this approach, see
\cite{Gio09}, section 4.3, or \cite{Gio10b}, section 12.

Another problem we have to take into account, and concerning roots
of infinitesimals, is that the equation $x^{2}=c$, for $c\in D_{\infty}$,
always has infinitely many solutions, e.g.

\begin{align*}
\left(\diff{t}_{4}\right)^{2} & =\diff{t}_{4/2}=\diff{t}_{2}\\
\left(\diff{t}_{4}+h\right)^{2} & =\diff{t}_{2}+h^{2}+2h\diff{t}_{4}=\diff{t}_{2}\quad\forall h\in D_{\infty}:\ \omega(h)<\frac{4}{3}.
\end{align*}

Therefore, we have infinitely many square roots of an infinitesimal. This
means that, although in $\R$ we have that $(-)^{2}:\R_{\ge0}\freccia\R_{\ge0}$
is bijective and $\sqrt{-}:\R_{\ge0}\freccia\R_{\ge0}$ is its (left
and right) inverse, in $\ER_{\ge0}$ we don't have injectivity, and,
therefore, we can have, at most, a right inverse. Indeed, we will
prove that $(-)^{2}:\ER_{\ge0}\freccia\ER_{\ge0}$ is surjective,
and
\[
\forall x\in\ER_{\ge0}:\ \left(\sqrt{x}\right)^{2}=x\quad\text{but}\quad\exists k:\ \sqrt{x^{2}}=_{k}x.
\]

Because we have infinitely many solutions of equations of the type $x^{p}=c$,
$p\in\R_{>1}$, a first idea is to choose, among them, the simplest
solution. Here, with {}``simplest'', we mean {}``the solution $x$
without unnecessary terms in its decomposition, i.e. without terms
that become zero taking the power $x^{\/p}$''. A similar idea of
{}``simplest solution'' has already been used in \cite{Gio09} to
define derivatives of smooth functions defined on infinitesimal sets.
For example, both $x=\diff{t}_{4}$ and $x=\diff{t}_{4}+\diff{t}$
are solutions of the equation $x^{2}=\diff{t}_{2}$, but, intuitively,
the first one is simpler compared to the second one, which contains
the unnecessary term $\diff{t}$.

However, there is another, more manageable idea to define roots of
infinitesimal numbers. Let
\[
c=\sum_{i=1}^{N_{c}}\st{c}_{i}\diff{t}_{\omega_{i}(c)}
\]
be the decomposition of $c\in\D_{\infty}$. Suppose $c>0$, so that
$\omega(c)=\omega_{1}(c)>\omega_{2}(c)>\dots>\omega_{N_{c}}(c)\ge1$
and $\st{c}_{i}\ne0$, $c_{1}>0$, then, for $p\in\R_{>0}$, we would
like to write

\begin{align}
c^{p} & =\left(\sum_{i=1}^{N_{c}}\st{c}_{i}\diff{t}_{\omega_{i}(c)}\right)^{p}=\nonumber \\
 & =\left[\st{c}_{1}\diff{t}_{\omega(c)}\cdot\left(1+\sum_{i=2}^{N_{c}}\frac{\st{c}_{i}}{\st{c}_{1}}\diff{t}_{\omega_{i}(c)\ominus\omega(c)}\right)\right]^{p}=\nonumber \\
 & =\left(\st{c}_{1}\right)^{p}\diff{t}_{\frac{\omega(c)}{p}}\cdot\left(1+\sum_{i=2}^{N_{c}}\frac{\st{c}_{i}}{\st{c}_{1}}\diff{t}_{\omega_{i}(c)\ominus\omega(c)}\right)^{p},\label{eq:idea_p_power}
\end{align}

where
\[
\frac{1}{0}:=\infty\quad,\quad a\oplus b:=\left(\frac{1}{a}+\frac{1}{b}\right)^{-1}\quad,\quad a\ominus b:=\left(\frac{1}{a}-\frac{1}{b}\right)^{-1}\quad\forall a,b\in\R.
\]
However, the right hand side of \eqref{eq:idea_p_power} is now a
well defined term, because the base of the $p$-th power is invertible.

\begin{rem}\ 

\begin{enumerate}
\item Note that the right hand side of \eqref{eq:idea_p_power} is well
defined if $\st{c}_{1}\ne0$, i.e. if $c\in D_{\infty}\setminus\{0\}$,
and because $\omega(c)>\omega_{i}(c)$, so that $\omega_{i}(c)\ominus\omega(c)>0$
and hence $\diff{t}_{\omega_{i}(c)\ominus\omega(c)}$ is well defined.
Moreover, it is not hard to prove that $\omega_{i}(c)\ominus\omega(c)>1$
if $i>1$.
\item It can be useful to note that setting
\begin{align*}
\ominus b & :=-b\\
a\odot b & :=a\cdot b,
\end{align*}
we easily have that $(\R_{\infty},\oplus,\ominus,\odot,\infty)$ is
a ring and the reciprocal function $\frac{1}{(-)}:\R_{\infty}\freccia\R$
is a ring isomorphism.
\end{enumerate}

\end{rem}

\begin{defn}
\label{def:c^p}Let $c\in D_{\infty}$, $c>0$, and $p\in\R_{>0}$,
then
\[
c^{p}:=\left(\st{c}_{1}\right)^{p}\diff{t}_{\frac{\omega(c)}{p}}\cdot\left(1+\sum_{i=2}^{N_{c}}\frac{\st{c}_{i}}{\st{c}_{1}}\diff{t}_{\omega_{i}(c)\ominus\omega(c)}\right)^{p}.
\]
Of course, if $p=\frac{m}{n}$, where $m$, $n\in\N$ and $n$ is
odd, the hypothesis $c>0$ can be dropped.
\end{defn}

\begin{example}\ 

\begin{enumerate}
\item Let us find $\sqrt{\diff{t}}$ using the previous definition. In this
case, we have $N_{c}=1$, $\st{c}_{1}=1$, $\omega(c)=1$, so
\[
\sqrt{\diff{t}}=(1)^{1/2}\diff{t}_{1/2}\cdot(1+0)^{1/2}=\diff{t}_{2}.
\]
\item We want to find $\sqrt{\diff{t}_{2}+\diff{t}}$: 

\begin{align*}
\sqrt{\diff{t}_{2}+\diff{t}} & =\diff{t}_{4}\cdot\left(1+\diff{t}_{1\ominus2}\right)^{1/2}=\diff{t}_{4}\left(1+\diff{t}_{2}\right)^{1/2}=\\
 & =\diff{t}_{4}\cdot\left(1+\sum_{n=1}^{2}\binom{\frac{1}{2}}{n}\diff{t}_{2/n}\right)=\\
 & =\diff{t}_{4}\cdot\left(1+\frac{1}{2}\cdot\diff{t}_{2}-\frac{1}{8}\cdot\diff{t}\right)=\\
 & =\diff{t}_{4}+\frac{1}{2}\diff{t}_{4/3}.
\end{align*}
We recall that

\begin{equation}
(1+h)^{p}=1+\sum_{n=1}^{+\infty}\binom{p}{n}\cdot h^{n}\quad\forall h\in D_{\infty}\label{eq:p_powerExpansion}
\end{equation}

because of the elementary transfer theorem (Theorem 24 in \cite{Gio10a}).
Finally, let us note that the series in \eqref{eq:p_powerExpansion}
is really a finite sum because of nilpotency of every infinitesimal
$h\in D_{\infty}$.\\
As expected, we have that $\left(\diff{t}_{4}+\frac{1}{2}\diff{t}_{4/3}\right)^{2}=\diff{t}_{2}+\diff{t}$.
\item $\sqrt{\diff{t}_{2}-\diff{t}_{3/2}-\diff{t}}=\diff{t}_{4}-\frac{1}{2}\diff{t}_{12/5}-\frac{1}{8}\diff{t}_{12/7}-\frac{9}{16}\diff{t}_{4/3}-\frac{37}{128}\diff{t}_{12/11}$.

\end{enumerate}

\end{example}

Generalizing these examples, we have that
\[
c^{p}=\left(\st{c}_{1}\right)^{p}\diff{t}_{\frac{\omega(c)}{p}}\cdot\left[1+\sum_{n=1}^{+\infty}\binom{p}{n}\cdot\left(\sum_{i=2}^{N_{c}}\frac{\st{c}_{i}}{\st{c}_{1}}\diff{t}_{\omega_{i}(c)\ominus\omega(c)}\right)^{n}\right].
\]

In the following theorems, in considering $x^{p}$ for generic $p\in\R_{>0}$,
we will always suppose $x>0$. However, this hypothesis can be dropped
in case of odd roots, and the proofs will remain essentially the same.

\begin{thm}
\label{thm:inverseUpToInf}Let $x\in D_{\infty}$, $x>0$, and $p\in\R$,
with $0<p<1$, then we have:

\begin{enumerate}
\item \label{enu:rightInverse}$\left(x^{p}\right)^{\frac{1}{p}}=x$
\item \label{enu:leftInverseUpToInf}If $x^{\frac{1}{p}}\ne0$ and $k:=\max\left\{ \omega_{2}(x),\omega_{2}\left[\left(x^{\frac{1}{p}}\right)^{p}\right]\right\} $,
then $\left(x^{\frac{1}{p}}\right)^{p}=_{k}x$.
\end{enumerate}

\end{thm}

\begin{rem}\ 

\begin{enumerate}
\item To understand better, it can be useful to clarify what is the difference
between the computation of $\left(x^{p}\right)^{1/p}$ and that of
$\left(x^{1/p}\right)^{p}$:

\begin{enumerate}
\item $\left(x^{p}\right)^{1/p}$: Because $p\in[0,1]_{\R}$, the computation
of $x^{p}$ is included in the Definition \ref{def:c^p}. Therefore,
we must:

\begin{enumerate}
\item Express $x$ using its decomposition.
\item Use Definition \ref{def:c^p}.
\item With the obtained result, we finally have to compute the subsequent
power $(-)^{1/p}$. However, $\frac{1}{p}>1$, so that this operation
is smooth and doesn't present any problem.
\end{enumerate}

\item $\left(x^{1/p}\right)^{p}$: In this case, the situation is the opposite
one.

\begin{enumerate}
\item $\frac{1}{p}>1$, so the operation $x^{1/p}=:y$ is smooth.
\item However, to compute $y^{p}$, we must apply Definition \ref{def:c^p},
so we firstly need the decomposition of $y=x^{1/p}$. Of course, it
is not easy to find this decomposition as a manageable function of
the decomposition of $x$.
\end{enumerate}

\end{enumerate}

\item Let us note that if $x$, $y\in D_{\infty}$ and $k:=\max\left[\omega(x),\omega(y)\right]$,
then it is trivially true that $x=_{k}y$, because $\omega(x-y)=k$
if $x\ne y$, and $\omega(x-y)=0$ otherwise. This shows that \emph{\ref{enu:leftInverseUpToInf}}
of Theorem \ref{thm:inverseUpToInf} is not trivial.
\item Theorem \ref{thm:inverseUpToInf} represents an overcoming of the
incompatibility between nil\-po\-tent infinitesimals and existence
of roots. Indeed, property \emph{\ref{enu:leftInverseUpToInf}} can
be applied only if $x^{1/p}\ne0$. Moreover, if $h\in D_{\infty}\setminus\{0\}$
and $h^{2}=0$, we have $\sqrt{h^{2}}=\sqrt{0}=0$, but, in general,
$\sqrt{h^{2}}\ne|h|$, e.g. $\sqrt{\left(\diff{t}\right)^{2}}=0\ne|\diff{t}|=\diff{t}$.

\end{enumerate}

\end{rem}

\noindent\textbf{Proof of Theorem \ref{thm:inverseUpToInf}:} Let $x=\sum_{j=1}^{N}b_{j}\diff{t}_{\beta_{j}}$
be the decomposition of $x$. Because $0<p<1$, from Definition \ref{def:c^p},
we have
\[
x^{p}=b_{1}^{p}\diff{t}_{\frac{\beta_{1}}{p}}\cdot\left(1+\sum_{j=2}^{N}\frac{b_{j}}{b_{1}}\diff{t}_{\beta_{j}\ominus\beta_{1}}\right)^{p}.
\]
Now, we have to apply the power $(-)^{1/p}$, which is smooth and
has the usual properties of powers (see, e.g., \cite{Gio09}, section
4.3). Therefore, we can write

\begin{align*}
\left(x^{p}\right)^{\frac{1}{p}} & =\left(b_{1}^{p}\diff{t}_{\frac{\beta_{1}}{p}}\right)^{\frac{1}{p}}\cdot\left[\left(1+\sum_{j=2}^{N}\frac{b_{j}}{b_{1}}\diff{t}_{\beta_{j}\ominus\beta_{1}}\right)^{p}\right]^{\frac{1}{p}}=\\
 & =b_{1}\diff{t}_{\beta_{1}}\cdot\left(1+\sum_{j=2}^{N}\frac{b_{j}}{b_{1}}\diff{t}_{\beta_{j}\ominus\beta_{1}}\right)=x.
\end{align*}
This proves \emph{\ref{enu:rightInverse}}.

To prove \emph{\ref{enu:leftInverseUpToInf},} we firstly have to
compute the smooth power $(-)^{1/p}$
\[
x^{\frac{1}{p}}=\left(\sum_{j=1}^{N}b_{j}\diff{t}_{\beta_{j}}\right)^{\frac{1}{p}}.
\]
The idea is to use the usual properties of $(-)^{1/p}$ and to gather
up the leading term $b_{1}\diff{t}_{\beta_{1}}$:

\begin{align}
x^{\frac{1}{p}} & =\left[b_{1}\diff{t}_{\beta_{1}}\cdot\left(1+\sum_{j=2}^{N}\frac{b_{j}}{b_{1}}\diff{t}_{\beta_{j}\ominus\beta_{1}}\right)\right]^{\frac{1}{p}}=\nonumber \\
 & =b_{1}^{\frac{1}{p}}\diff{t}_{p\beta_{1}}\cdot\left(1+\sum_{j=2}^{N}\frac{b_{j}}{b_{1}}\diff{t}_{\beta_{j}\ominus\beta_{1}}\right)^{\frac{1}{p}}.\label{eq:gatherLeadingTerm}
\end{align}

We are not able to find the decomposition of this number, but we can
surely claim that

\begin{align}
\st{\left(x^{\frac{1}{p}}\right)}_{1} & =b_{1}^{\frac{1}{p}}=\left(\st{x}_{1}\right)^{\frac{1}{p}}\label{eq:2_firstStdPart}\\
\omega\left(x^{\frac{1}{p}}\right) & =p\cdot\beta_{1}=p\cdot\omega(x)\quad\text{if }x^{\frac{1}{p}}\ne0\text{ i.e. if }p\cdot\beta_{1}\ge1.\label{eq:2_order}
\end{align}

This guess is based on the idea that in \eqref{eq:gatherLeadingTerm},
the infinitesimal $b_{1}^{\frac{1}{p}}\cdot\diff{t}_{p\beta_{1}}$
is multiplied by an invertible number, whose standard part is 1. Indeed,
we have the following

\begin{lem}
\label{lem:orderAndStdPartOfProductByInvertible}Let $k$, $h\in D_{\infty}$
and $y\in\ER$ such that $y$ is invertible and $k=h\cdot y$. Then

\begin{align*}
\omega(k) & =\omega(h)\\
\st{k}_{1} & =\st{h}_{1}\cdot\st{y}.
\end{align*}

\end{lem}

We postpone the proof of this lemma to the end of the current proof.

Using \eqref{eq:2_firstStdPart} and \eqref{eq:2_order} and Definition
\ref{def:c^p}, we can write

\begin{equation}
\left(x^{\frac{1}{p}}\right)^{p}=\st{x}_{1}\diff{t}_{\omega(x)}\cdot\left(1+\sum_{i=2}^{M}\frac{\st{\left(x^{\frac{1}{p}}\right)}_{i}}{\st{x}_{1}}\diff{t}_{\omega_{i}\left(x^{\frac{1}{p}}\right)\ominus\omega(x)}\right)^{p},\label{eq:3_x_1overp_p}
\end{equation}
where $M$ is the number of terms in the decomposition of $x^{1/p}$
(of which, we really know only the first term). Using again Lemma
\ref{lem:orderAndStdPartOfProductByInvertible} applied to \eqref{eq:3_x_1overp_p},
we have

\begin{align*}
\st{\left[\left(x^{\frac{1}{p}}\right)^{p}\right]}_{1} & =\st{x}_{1}\\
\omega\left[\left(x^{\frac{1}{p}}\right)^{p}\right] & =\omega(x).
\end{align*}

Let us observe that $x>0$, hence $x\ne0$ and $\omega(x)\ge1$, so
that $\omega\left(\st{x}_{1}\diff{t}_{\omega(x)}\right)=\omega(x)$.
We have hence proved that the first terms in the decompositions of
both $x$ and $\left(x^{\frac{1}{p}}\right)^{p}$ are the same. Therefore
\[
\omega\left[x-\left(x^{\frac{1}{p}}\right)^{p}\right]\le\max\left\{ \omega_{2}(x),\omega_{2}\left[\left(x^{\frac{1}{p}}\right)^{p}\right]\right\} =k
\]
 and hence $x=_{k}\left(x^{\frac{1}{p}}\right)^{p}$.$\qed$

\noindent \textbf{Proof of Lemma \ref{lem:orderAndStdPartOfProductByInvertible}:}
Because $y$ is invertible, we have that $\st{y}\ne0$. Write the
product $k=h\cdot y$ using decompositions

\begin{align}
k & =\sum_{i=1}^{N_{k}}\st{k}_{i}\diff{t}_{\omega_{i}(k)}=\left(\sum_{p=1}^{N_{h}}\st{h}_{p}\diff{t}_{\omega_{p}(h)}\right)\cdot\left(\st{y}+\sum_{q=1}^{N_{y}}\st{y}_{q}\diff{t}_{\omega_{q}(y)}\right)=\nonumber \\
 & =\st{h}_{1}\st{y}\diff{t}_{\omega(h)}+\sum_{p=2}^{N_{h}}\st{h}_{p}\st{y}\diff{t}_{\omega_{p}(h)}+\sum_{p=1}^{N_{h}}\sum_{q=1}^{N_{y}}\st{h}_{p}\st{y}_{q}\diff{t}_{\omega_{p}(h)\oplus\omega_{q}(y)}.\label{eq:productByDecomp}
\end{align}

However, in general, we have $a\oplus b<\min(a,b)$ for every $a$,
$b\in\R_{>0}$, so that, in \eqref{eq:productByDecomp} the leading
term is $\st{h}_{1}\st{y}\diff{t}_{\omega(h)}$ and hence from the
uniqueness of decomposition, the conclusion follows.$\qed$

\begin{rem}\ 

\begin{enumerate}
\item Let us observe that, in the hypothesis of Theorem \ref{thm:inverseUpToInf},
we also have
\[
y=x^{p}\then\left(y^{\frac{1}{p}}\right)^{p}=y.
\]
In fact, $y^{\frac{1}{p}}=\left(x^{p}\right)^{\frac{1}{p}}=x$, and
hence $\left(y^{\frac{1}{p}}\right)^{p}=x^{p}=y$. Therefore, for
numbers of the form $y=x^{p}$, the equality \emph{\ref{enu:leftInverseUpToInf}}
of Theorem \ref{thm:inverseUpToInf} becomes exact. We can interpret
this result saying that our Definition \ref{def:c^p} of $c^{p}$
gives exactly the simplest solution of the equation $x^{\frac{1}{p}}=c$.
Indeed, like in the case $\sqrt{\left(\diff{t}_{2}+\diff{t}\right)^{2}}=_{1}\diff{t}_{2}+\diff{t}$,
we can say that in \emph{\ref{enu:leftInverseUpToInf}} we don't have
an exact equality if the number $x$ contains unnecessary infinitesimals
with respect to the power $(-)^{1/p}$, like $\diff{t}$ in the previous
example. See section \ref{sub:The-notion-of-incomplete-term} for
a formalization of the notion of {}``unnecessary term with respect
to the power $(-)^{1/p}$''.
\item The equality \emph{\ref{enu:leftInverseUpToInf},} up to infinitesimals,
implies that
\[
\sqrt{-}:\ER_{\ge0}\freccia\ER_{\ge0}\quad\text{is not surjective.}
\]
For example $y=\diff{t}_{2}+\diff{t}$ cannot be written as the square
root of some number $x$. Otherwise, we would have
\[
y=\sqrt{x}\then y^{2}=\left(\diff{t}_{2}+\diff{t}\right)^{2}=\diff{t}=\left(\sqrt{x}\right)^{2}=x,
\]
but then $\sqrt{x}=\sqrt{\diff{t}}=\diff{t}_{2}\ne y$. Of course,
this corresponds to saying that the square is not the right inverse of
the square root.
\item Trivially, we can consider a smooth function $g:\R\freccia\R$ having
a root of order $n\in\N_{>0}$ at $x=0$, i.e. such that
\[
g(h)=a\cdot h^{n}\quad\forall h\in D_{n},
\]
where $a\in\R_{\ne0}$. We can hence define a sort of infinitesimal
right inverse of $g$, setting
\[
f(h):=\sqrt[n]{\frac{h}{a}}\quad\forall h\in D_{\infty}
\]
if $n$ is odd or $n$ is even and $a>0$, and
\[
f(h):=\sqrt[n]{\frac{h}{-a}}\quad\forall h\in D_{\infty}
\]
if $n$ is even and $a<0$. Then we have $g(f(h))=\pm h$ for every
$h\in D_{\infty}$, with the positive sign in the first case.

\end{enumerate}

\end{rem}

\subsection{A formula to compute a root}

By definition, if $c\in D_{\infty}\setminus\{0\}$, we have
\[
c^{p}:=\left(\st{c}_{1}\right)^{p}\diff{t}_{\frac{\omega(c)}{p}}\cdot\left(1+\sum_{i=2}^{N_{c}}\frac{\st{c}_{i}}{\st{c}_{1}}\diff{t}_{\omega_{i}(c)\ominus\omega(c)}\right)^{p}.
\]
The $p$-th power of the invertible term can be computed in several,
obviously equivalent, ways. 

\begin{enumerate}
\item Using the infinitesimal Taylor formula of the function $(1+x)^{p}$,
with

\begin{equation}
x=\sum_{i=2}^{N_{c}}\frac{\st{c}_{i}}{\st{c}_{1}}\diff{t}_{\omega_{i}(c)\ominus\omega(c)}\in D_{\omega_{2}(c)\ominus\omega_{1}(c)}.\label{eq:infTermInFormula}
\end{equation}

\item Equivalently, we can use the formula $(1+x)^{p}=1+\sum_{n=1}^{+\infty}\binom{p}{n}\cdot x^{n}$,
for $|x|<1$, which transfers to $D_{\infty}$ by Theorem 24 of \cite{Gio10a}.

\end{enumerate}

Applying the second method, we get

\begin{align*}
c^{p} & =\left(\st{c}_{1}\right)^{p}\diff{t}_{\frac{\omega(c)}{p}}\cdot\left[1+\sum_{n=1}^{+\infty}\binom{p}{n}\cdot\left(\sum_{i=2}^{N_{c}}\frac{\st{c}_{i}}{\st{c}_{1}}\diff{t}_{\omega_{i}(c)\ominus\omega(c)}\right)^{n}\right]=\\
 & =\left(\st{c}_{1}\right)^{p}\diff{t}_{\frac{\omega(c)}{p}}\cdot\left[1+\sum_{n=1}^{k_{c,p}}\binom{p}{n}\sum_{\substack{\gamma\in\N^{N_{c}-1}\\
|\gamma|=n}
}\frac{n!}{\gamma!}\cdot\prod_{i=2}^{N_{c}}\left(\frac{\st{c}_{i}}{\st{c}_{1}}\right)^{\gamma_{i-1}}\diff{t}_{\frac{\omega_{i}(c)\ominus\omega(c)}{\gamma_{i-1}}}\right],
\end{align*}
where

\begin{equation}
k_{c,p}:=
\begin{cases}
\ulcorner\omega_{2}(c)-\omega(c)\urcorner & \text{ if }p\notin\N\\
\min\left(\ulcorner\omega_{2}(c)-\omega(c)\urcorner,p\right) & \text{ if }p\in\N
\end{cases}\label{eq:defOf_kcp}
\end{equation}
with $\ulcorner a\urcorner\in\N$ the ceiling of $a\in\R$, that
is the smallest integer greater than or equal to $a$. Note that the
first alternative of \eqref{eq:defOf_kcp} is due to \eqref{eq:infTermInFormula},
whereas the second one is also due to the equality $\binom{p}{n}=0$
if $n>p$.

Using Theorem 13 of \cite{Gio10a}, we have

\begin{align*}
\prod_{i=2}^{N_{c}}\diff{t}_{\frac{\omega_{i}(c)\ominus\omega(c)}{\gamma_{i-1}}}\ne0 & \iff\sum_{i=2}^{N_{c}}\frac{\gamma_{i-1}}{\omega_{i}(c)\ominus\omega(c)}\le1\\
 & \iff\bigoplus_{i=2}^{N_{c}}\frac{\omega_{i}(c)\ominus\omega(c)}{\gamma_{i-1}}\ge1.
\end{align*}

In the following, we will set $\omega(c,\gamma):=\bigoplus_{i=2}^{N_{c}}\frac{\omega_{i}(c)\ominus\omega(c)}{\gamma_{i-1}}$,
so that
\[
\prod_{i=2}^{N_{c}}\diff{t}_{\frac{\omega_{i}(c)\ominus\omega(c)}{\gamma_{i-1}}}=\diff{t}_{\omega(c,\gamma)}.
\]
Finally, we obtain the formula

\begin{equation}
c^{p}=\left(\st{c}_{1}\right)^{p}\diff{t}_{\frac{\omega(c)}{p}}+\sum_{n=1}^{k_{c,p}}\binom{p}{n}\sum_{\gamma\in\Gamma_{c,n}}\frac{n!}{\gamma!}\cdot\st{c}_{1}^{p-n}\cdot\st{c}_{2}^{\gamma_{1}}\cdot\ldots\cdot\st{c}_{N_{c}}^{\gamma_{N_{c}-1}}\diff{t}_{\omega(c,\gamma)\oplus\frac{\omega(c)}{p}},\label{eq:finalFormula}
\end{equation}
where
\[
\Gamma_{c,n}:=\left\{ \gamma\in\N^{N_{c}-1}\,|\,|\gamma|=n\ ,\ \omega(c,\gamma)\ge1\right\} .
\]

\subsection{Properties of roots: the general theorem}

For generic $x$, $y\in\ER$ we can state the following

\begin{thm}
\label{thm:propertiesOfPowersGeneric}Let $x$, $y\in\D_{\infty}$
be strictly positive infinitesimals, and $p$, $q\in\R_{>0}$, then:

\begin{enumerate}
\item \label{enu:order1}$\omega\left[\left(x^{p}\right)^{q}\right]=\omega(x^{pq})=:o_{1}$
and $\st{\left[\left(x^{p}\right)^{q}\right]}_{1}=\st{\left(x^{pq}\right)}_{1}$.
\item \label{enu:order2}$\omega\left[\left(x\cdot y\right)^{p}\right]=\omega(x^{p}\cdot y^{p})=:o_{2}$
and $\st{\left[\left(x\cdot y\right)^{p}\right]}_{1}=\st{\left(x^{p}\cdot y^{p}\right)}_{1}$.
\item \label{enu:prop1UpToInf}$\exists k\in\R:\ 1\le k<o_{1}$ and $\left(x^{p}\right)^{q}=_{k}x^{pq}$.
\item \label{enu:prop2UpToInf}$\exists k\in\R:\ 1\le k<o_{2}$ and $\left(x\cdot y\right)^{p}=_{k}x^{p}\cdot y^{p}$.
\item \label{enu:prop3}$x^{p}\cdot y^{q}=x^{p+q}$.
\end{enumerate}

\end{thm}

Before proving this theorem, we need the following very useful lemma:

\begin{lem}
\label{lem:formulaForPowersWithoutUsingDecomp}Let $c=\sum_{j=1}^{M}a_{j}\diff{t}_{\alpha_{j}}$,
with $\alpha_{1}>\alpha_{j}\ge1$ for every $j=2,\dots,M$, and $a_{1}>0$.
Let us note explicitly that not necessarily this is the decomposition
of $c$. Then
\[
c^{p}=a_{1}^{p}\diff{t}_{\frac{\alpha_{1}}{p}}\cdot\left(1+\sum_{j=2}^{M}\frac{a_{j}}{a_{1}}\diff{t}_{\alpha_{j}\ominus\alpha_{1}}\right)^{p}.
\]
\end{lem}

Of course, this lemma states that the formula used for the definition
of $c^{p}$ can also be used starting from a representation $c=\sum_{j=1}^{M}a_{j}\diff{t}_{\alpha_{j}}$
which is not necessarily the decomposition of $c$. To apply this
lemma, the important step is to find the greatest infinitesimal $a_{1}\diff{t}_{\alpha_{1}}$
and to check that all the other terms $a_{j}\diff{t}_{\alpha_{j}}$
can be zero only if $a_{j}=0$.

\noindent \textbf{Proof of Lemma \ref{lem:formulaForPowersWithoutUsingDecomp}:}
Starting from $c=\sum_{j=1}^{M}a_{j}\diff{t}_{\alpha_{j}}$, we firstly
sum all the coefficients $a_{j}$ having the same infinitesimal $\diff{t}_{\alpha_{j}}$,
i.e. if

\begin{align}
\bar{a}_{q} & :=\sum\left\{ a_{j}\,|\, j=1,\dots,M\ ,\ \alpha_{j}=q\right\} \quad\forall q\in\left\{ \alpha_{j}\,|\, j=1,\dots,M\right\} \nonumber \\
O & :=\left\{ \alpha_{j}\,|\, j=1,\dots,M\ ,\ \bar{a}_{\alpha_{j}}\ne0\right\} =:\left\{ q_{1},\dots,q_{N}\right\} ,\label{eq:enumerationOfOrders}
\end{align}

then

\begin{equation}
c=\sum_{q\in O}\bar{a}_{q}\diff{t}_{q}=\sum_{i=1}^{N}\bar{a}_{q_{i}}\diff{t}_{q_{i}}.\label{eq:c_sumOfAlpha}
\end{equation}

Let us note that in \eqref{eq:enumerationOfOrders}, $q_{1},\dots,q_{N}$
is any enumeration of the elements of the set of all orders $O$.
Now, all the summands in \eqref{eq:c_sumOfAlpha} are non zero, because
of our definition of the set $O$. Therefore, reordering the summands
in \eqref{eq:c_sumOfAlpha}, we obtain the decomposition of $c$.
Formally, this means that we can find a permutation $\sigma$ of $\{1,\dots,N\}$
such that

\begin{equation}
c=\sum_{i=1}^{N}\bar{a}_{q_{\sigma_{i}}}\diff{t}_{q_{\sigma_{i}}}\label{eq:decompositionOf_c}
\end{equation}

is the decomposition of $c$. Let us note that, to obtain \eqref{eq:decompositionOf_c},
we need that for every $i=1,\dots N$ we can find $j=1,\dots M$ such
that $q_{i}=\alpha_{j}\ge1$. By definition of decomposition, $q_{\sigma_{1}}$
is the maximum order in \eqref{eq:decompositionOf_c}, i.e. $q_{\sigma_{1}}=\max\left\{ q_{1},\dots,q_{N}\right\} =\max\left\{ \alpha_{j}\,|\,\alpha_{j}\in O\right\} $.
However, we have that $q_{\sigma_{1}}=\alpha_{1}$, because $\alpha_{1}>\alpha_{j}$
by hypothesis, and because
\[
\bar{a}_{\alpha_{1}}=\sum\left\{ a_{j}\,|\, j=1,\dots,M\ ,\ \alpha_{j}=\alpha_{1}\right\} =a_{1}\ne0,
\]
so that $\alpha_{1}\in O$ and so $\bar{a}_{q_{\sigma_{1}}}=\bar{a}_{\alpha_{1}}=a_{1}$.
We can now apply our Definition \ref{def:c^p} using the decomposition
\eqref{eq:decompositionOf_c}:
\[
c^{p}=a_{1}^{p}\diff{t}_{\frac{\alpha_{1}}{p}}\cdot\left(1+\sum_{i=2}^{N}\frac{\bar{a}_{q_{\sigma_{i}}}}{a_{1}}\diff{t}_{q_{\sigma_{i}}\ominus\alpha_{1}}\right)^{p}.
\]
Now, we only have to retrace the previous steps, so as to eliminate
$\sigma$, $q$, $\bar{a}$, etc.

\begin{align*}
c^{p} & =a_{1}^{p}\diff{t}_{\frac{\alpha_{1}}{p}}\cdot\left(1+\sum_{i=2}^{N}\frac{\bar{a}_{q_{i}}}{a_{1}}\diff{t}_{q_{i}\ominus\alpha_{1}}\right)^{p}=\\
 & =a_{1}^{p}\diff{t}_{\frac{\alpha_{1}}{p}}\cdot\left(1+\sum_{\substack{q\in O\\
q\ne\alpha_{1}}
}\frac{\bar{a}_{q}}{a_{1}}\diff{t}_{q\ominus\alpha_{1}}\right)^{p}=\\
 & =a_{1}^{p}\diff{t}_{\frac{\alpha_{1}}{p}}\cdot\left(1+\sum_{\substack{j=2\\
\bar{a}_{\alpha_{j}}\ne0}
}^{M}\frac{a_{j}}{a_{1}}\diff{t}_{\alpha_{j}\ominus\alpha_{1}}+\sum_{\substack{j=2\\
\bar{a}_{\alpha_{j}}=0}
}^{M}\frac{0}{a_{1}}\diff{t}_{\alpha_{j}\ominus\alpha_{1}}\right)^{p}=\\
 & =a_{1}^{p}\diff{t}_{\frac{\alpha_{1}}{p}}\cdot\left(1+\sum_{i=2}^{N}\frac{a_{j}}{a_{1}}\diff{t}_{\alpha_{j}\ominus\alpha_{1}}\right)^{p},
\end{align*}

which is our conclusion.$\qed$

\noindent \textbf{Proof of Theorem \ref{thm:propertiesOfPowersGeneric}:}
To prove \emph{\ref{enu:order1},} let $x=\sum_{j=1}^{N}b_{j}\diff{t}_{\beta_{j}}$
be the decomposition of $x$. The idea is to use formula \eqref{eq:finalFormula}
to compute $x^{p}$, and then Lemma \ref{lem:formulaForPowersWithoutUsingDecomp}
to compute $\left(x^{p}\right)^{q}$. To avoid heavy notations, we
will use the simplified symbols $k:=k_{x,p}$ and $\Gamma:=\Gamma_{x,n}$:

\begin{align*}
\left(x^{p}\right)^{q} & =\left[b_{1}^{p}\diff{t}_{\frac{\beta_{1}}{p}}+\right.\\
 & \phantom{=}\left.+\sum_{n=1}^{k}\binom{p}{n}\sum_{\gamma\in\Gamma}\frac{n!}{\gamma!}\cdot b_{1}^{p-1}\cdot b_{2}^{\gamma_{1}}\cdot\ldots\cdot b_{N_{c}}^{\gamma_{N_{c}-1}}\diff{t}_{\omega(x,\gamma)\oplus\frac{\beta_{1}}{p}}\right]^{q}.
\end{align*}

Before using Lemma \ref{lem:formulaForPowersWithoutUsingDecomp},
we need to prove that $\diff{t}_{\frac{\beta_{1}}{p}}$ is the greatest
infinitesimal, so let us compute

\begin{align*}
\left[\omega(x,\gamma)\oplus\frac{\beta_{1}}{p}\right]^{-1} & =\frac{p}{\beta_{1}}+\sum_{i=2}^{N}\frac{\gamma_{i-1}}{\beta_{i}\ominus\beta_{1}}=\\
 & =\frac{p}{\beta_{1}}+\sum_{i=2}^{N}\gamma_{i-1}\left(\frac{1}{\beta_{i}}-\frac{1}{\beta_{1}}\right)=\\
 & =\frac{p-n}{\beta_{1}}+\sum_{i=2}^{N}\frac{\gamma_{i-1}}{\beta_{i}}.
\end{align*}

So, we need to prove that $\frac{p}{\beta_{1}}<\frac{p-n}{\beta_{1}}+\sum_{i=2}^{N}\frac{\gamma_{i-1}}{\beta_{i}}$,
that is $n<\sum_{i=2}^{N}\gamma_{i-1}\cdot\frac{\beta_{1}}{\beta_{i}}$.
In fact, $\beta_{1}>\beta_{i}$ so that $\sum_{i=2}^{N}\gamma_{i-1}\cdot\frac{\beta_{1}}{\beta_{i}}>\sum_{i=2}^{N}\gamma_{i-1}=n$.
Moreover, we can suppose to restrict the set $\Gamma$ to those $\gamma$
such that $\omega(x,\gamma)\oplus\frac{\beta_{1}}{p}\ge1$, because,
otherwise, the corresponding term $\diff{t}_{\omega(x,\gamma)\oplus\frac{\beta_{1}}{p}}=0$.
We can hence apply the Lemma \ref{lem:formulaForPowersWithoutUsingDecomp},
obtaining

\begin{multline}
\left(x^{p}\right)^{q}=b_{1}^{pq}\diff{t}_{\frac{\beta_{1}}{pq}}\cdot\\
\cdot\left[1+\sum_{n=1}^{k}\binom{p}{n}\sum_{\gamma\in\Gamma}\frac{n!}{\gamma!}\cdot b_{1}^{p-n}\cdot b_{2}^{\gamma_{1}}\cdot\ldots\cdot b_{N_{c}}^{\gamma_{N_{c}-1}}\diff{t}_{\omega(x,\gamma)\oplus\frac{\beta_{1}}{p}}\right]^{q}.\label{eq:1_xpq}
\end{multline}
On the other hand, we have

\begin{equation}
x^{pq}=b_{1}^{pq}\diff{t}_{\frac{\beta_{1}}{pq}}\cdot\left(1+\sum_{j=2}^{N}\frac{b_{j}}{b_{1}}\diff{t}_{\beta_{j}\ominus\beta_{1}}\right)^{pq}.\label{eq:2_xpq}
\end{equation}
Therefore, the conclusion follows from \eqref{eq:1_xpq}, \eqref{eq:2_xpq}
and Lemma \ref{lem:orderAndStdPartOfProductByInvertible}.

To prove \emph{\ref{enu:order2}}, we can use the same method as
before. Let $y=\sum_{k=1}^{M}e_{k}\diff{t}_{\varepsilon_{k}}$ be the
decomposition of $y$, then
\[
\left(xy\right)^{p}=\left(\sum_{j=1}^{N}\sum_{k=1}^{M}b_{j}e_{k}\diff{t}_{\beta_{j}\oplus\eps_{k}}\right)^{p}.
\]
Of course, $\beta_{1}\oplus\eps_{1}>\beta_{j}\oplus\eps_{k}$ for
any $j$ and $k$, and we can also consider
\[
I:=\left\{ (j,k)\,|\, j=1,\dots,N\ ,\ k=1,\dots,M\ ,\ \beta_{j}\oplus\eps_{k}\ge1\right\} ,
\]
so that to the sum
\[
\left(xy\right)^{p}=\left(\sum_{(j,k)\in I}b_{j}e_{k}\diff{t}_{\beta_{j}\oplus\eps_{k}}\right)^{p}
\]
we can apply Lemma \ref{lem:formulaForPowersWithoutUsingDecomp}.
We obtain
\[
\left(xy\right)^{p}=b_{1}^{p}e_{1}^{p}\diff{t}_{\frac{\beta_{1}\oplus\eps_{1}}{p}}\cdot\left(1+\sum_{\substack{(j,k)\in I\\
(j,k)\ne(1,1)}
}\frac{b_{j}\cdot e_{k}}{b_{1}\cdot e_{1}}\diff{t}_{\beta_{j}\oplus\eps_{k}\ominus\left(\beta_{1}\oplus\eps_{1}\right)}\right)^{p}.
\]
On the other hand, we have

\begin{multline*}
x^{p}\cdot y^{p}=b_{1}^{p}\diff{t}_{\frac{\beta_{1}}{p}}\cdot\left(1+\sum_{j=2}^{N}\frac{b_{j}}{b_{1}}\diff{t}_{\beta_{j}\ominus\beta_{1}}\right)^{p}\cdot\\
\cdot e_{1}^{p}\diff{t}_{\frac{\eps_{1}}{p}}\cdot\left(1+\sum_{k=2}^{M}\frac{e_{k}}{e_{1}}\diff{t}_{\eps_{k}\ominus\eps_{1}}\right)^{p}=b_{1}^{p}e_{1}^{p}\diff{t}_{\frac{\beta_{1}\oplus\eps_{1}}{p}}\cdot\left(1+h\right)^{p},
\end{multline*}
where $h\in D_{\infty}$ is obtained from the product of the previous
$p$-th powers with invertible bases. Once again, the conclusion follows
from Lemma \ref{lem:orderAndStdPartOfProductByInvertible}.

The proofs of \emph{\ref{enu:prop1UpToInf}} and \emph{\ref{enu:prop2UpToInf}}
are straightforward, taking into account \emph{\ref{enu:order1}} and \emph{\ref{enu:order2}}
so that, e.g., in the difference $\left(x^{p}\right)^{q}-x^{pq}$
there appear only infinitesimals of order greater than $o_{1}$.

Finally, to prove \emph{\ref{enu:prop3}}, we only have to apply our
Definition \ref{def:c^p} of power:

\begin{align*}
x^{p}\cdot x^{q} & =b_{1}^{p}\diff{t}_{\frac{\beta_{1}}{p}}\cdot\left(1+\sum_{j=2}^{N}\frac{b_{j}}{b_{1}}\diff{t}_{\beta_{j}\ominus\beta_{1}}\right)^{p}\cdot\\
 & \phantom{=}\cdot b_{1}^{q}\diff{t}_{\frac{\beta_{1}}{q}}\cdot\left(1+\sum_{j=2}^{N}\frac{b_{j}}{b_{1}}\diff{t}_{\beta_{j}\ominus\beta_{1}}\right)^{q}=\\
 & =b_{1}^{p+q}\diff{t}_{\frac{\beta_{1}}{p+q}}\cdot\left(1+\sum_{j=2}^{N}\frac{b_{j}}{b_{1}}\diff{t}_{\beta_{j}\ominus\beta_{1}}\right)^{p+q}=\\
 & =x^{p+q},
\end{align*}

which is the conclusion.$\qed$

\begin{rem}
For generic $x$, $y\in D_{\infty}$, properties \emph{\ref{enu:prop1UpToInf}}
and \emph{\ref{enu:prop2UpToInf}} of Theorem \ref{thm:propertiesOfPowersGeneric}
cannot be improved. Indeed, if we had always
\[
\left(x^{p}\right)^{q}=x^{pq}\quad\forall x\in D_{\infty}\ \forall p,q>0
\]
we would have, as a consequence, the general validity of 
\[
\left(x^{\frac{1}{p}}\right)^{p}=x^{\frac{1}{p}p}=x\quad\forall x\in D_{\infty}\ \forall p\in(0,1]_{\R},
\]
but we know that this property is not generally true.

Analogously, from the general validity of
\[
(xy)^{p}=x^{p}\cdot y^{p}\quad\forall x\in D_{\infty}\ \forall p>0,
\]
we would have
\[
\left(x^{2}\right)^{\frac{1}{2}}=(x\cdot x)^{\frac{1}{2}}=x^{\frac{1}{2}}\cdot x^{\frac{1}{2}}=\left(\sqrt{x}\right)^{2}=x\quad\forall x\in D_{\infty},
\]
but we know that this is not generally true.

\end{rem}

Therefore, on the one hand these seem the best results attainable.
However, it does not seem desirable to work with the
equality $=_{k}$ up to infinitesimals of some order, in particular
for such basic operations.

We have already noted that several counterexamples are of the form
\[
\left[\left(\diff{t}_{2}+\diff{t}\right)^{2}\right]^{\frac{1}{2}}=\diff{t}_{2},
\]
where we have terms like $\diff{t}$ which are, intuitively, unnecessary
with respect to the square. Our next aim is to formalize the idea
of \emph{incomplete term with respect to $(-)^{\frac{1}{p}}$}, and
to prove that the usual properties of the powers hold, with the usual
equality, if we use only Fermat reals without incomplete terms. We
will also see why the name \emph{incomplete term} seems a better choice
than \emph{unnecessary term}.

\subsection{The notion of incomplete term\label{sub:The-notion-of-incomplete-term}}

Let us start from the usual notations and hypotheses: $x=\sum_{j=1}^{N}b_{j}\diff{t}_{\beta_{j}}$
is the decomposition of $x$, and $p\in(0,1]_{\R}$. In this decomposition,
let us consider a term $\diff{t}_{\beta_{r+1}}$, for $r=1,\dots,N-1$.
The power $(-)^{1/p}$ is smooth, because $\frac{1}{p}>1$, and, with
the usual calculations, we can write
\begin{equation}
x^{\frac{1}{p}}=b_{1}^{\frac{1}{p}}\diff{t}_{p\beta_{1}}+\sum_{n=1}^{k}\binom{\frac{1}{p}}{n}\sum_{\substack{\gamma\in\N^{N-1}\\
|\gamma|=n}
}\frac{n!}{\gamma!}\cdot b_{1}^{\frac{1}{p}-n}\cdot b_{2}^{\gamma_{1}}\cdot\ldots\cdot b_{N}^{\gamma_{N-1}}\diff{t}_{\omega(x,\gamma)\oplus p\beta_{1}},\label{eq:x_1over_p_powerFormula}
\end{equation}
where
\begin{align*}
k:=k_{x,\frac{1}{p}} & :=
\begin{cases}
\ulcorner\beta_{2}\ominus\beta_{1}\urcorner & \text{ if }\frac{1}{p}\notin\N\\
\min\left(\ulcorner\beta_{2}\ominus\beta_{1}\urcorner,\frac{1}{p}\right) & \text{ if }\frac{1}{p}\in\N
\end{cases}\\
\omega(x,\gamma) & :=\bigoplus_{j=2}^{N}\frac{\beta_{j}\ominus\beta_{1}}{\gamma_{i-1}}.
\end{align*}

We have two possibilities to identify the terms, like $\diff{t}$
in $\left(\diff{t}_{2}+\diff{t}\right)^{2}$, that are \emph{unnecessary},
or, better, \emph{incomplete}.

The first one is to say that a term of the type $\diff{t}_{\beta_{r+1}}$
\emph{gives no contribution} whenever expanding the power $x^{1/p}$,
it gives \emph{always} zero summands, exactly like $\diff{t}$ in
$\left(\diff{t}_{2}+\diff{t}\right)^{2}=\diff{t}+\left(\diff{t}\right)^{2}+2\cdot\diff{t}_{2}\cdot\diff{t}=\diff{t}$.
Putting it in negative form: if, expanding the power $x^{1/p}$, we
have that \emph{at least one summand}, obtained from $\diff{t}_{\beta_{r+1}}$,
is not zero, then the term $\diff{t}_{\beta_{r+1}}$ \emph{gives}
\emph{some} \emph{contribution}, i.e. it is necessary.

A substantial objection against this idea, however, is the following: let us
suppose that $\diff{t}_{\beta_{r+1}}$ is necessary, i.e. it gives some contribution.
Then, the situation described above also includes the possibility
that, in the expansion of the power $x^{1/p}$, the term $\diff{t}_{\beta_{r+1}}$
gives, e.g., only one contribution, whereas all the other terms involving
$\diff{t}_{\beta_{r+1}}$ give zero. One of our first aims will be
to prove that, if in the decomposition of $x$ every term gives a
contribution, then $\left(x^{1/p}\right)^{p}=x$. The situation can
actually be problematic because, following the previous extreme example,
in $x^{1/p}$ {}``all the information concerning $\diff{t}_{\beta_{r+1}}$''
is contained in the unique non zero term. For example,
\begin{align*}
x & =\diff{t}_{3}+\diff{t}_{\frac{3}{2}}\quad,\quad\diff{t}_{\beta_{r+1}}=\diff{t}_{\frac{3}{2}}\\
x^{2} & =\diff{t}_{\frac{3}{2}}+\diff{t}_{\frac{3}{4}}+2\diff{t}_{3\oplus\frac{3}{2}}=\diff{t}_{\frac{3}{2}}+2\diff{t}.
\end{align*}
Can the inverse operation $x^{2}\mapsto\sqrt{x^{2}}$ reconstruct
the whole initial information, about $x$, starting only from the
unique non zero term $2\diff{t}$ generated by $\diff{t}_{\beta_{r+1}}$?

\begin{align*}
\sqrt{x^{2}} & =\diff{t}_{3}\cdot\left(1+2\diff{t}_{1\ominus\frac{3}{2}}\right)^{\frac{1}{2}}=\\
 & =\diff{t}_{3}\cdot\left(1+\sum_{n=1}^{+\infty}\binom{\frac{1}{2}}{n}\cdot2^{n}\diff{t}_{\frac{3}{n}}\right)=\\
 & =\diff{t}_{3}\cdot\left(1+\frac{1}{2}2\diff{t}_{3}-\frac{1}{8}4\diff{t}_{\frac{3}{2}}+\frac{1}{16}8\diff{t}\right)=\\
 & =\diff{t}_{3}+\diff{t}_{\frac{3}{2}}-\frac{1}{2}\diff{t}.
\end{align*}
This counterexample hence gives a negative answer to our question. 

The second possibility of defining a precise notion of incomplete term
arises from trying to prove the property $\left(x^{1/p}\right)^{p}=x$
starting from a definition based on the previous erroneous idea. We
will say, intuitively, that $\diff{t}_{\beta_{r+1}}$ is incomplete
whenever expanding the power $x^{1/p}$ gives \emph{at least
one} zero summand. For this reason, the term {}``incomplete'' is
better than {}``unnecessary''. Putting it in negative form: if expanding
the power $x^{1/p}$, we have that \emph{every summand}, obtained
from $\diff{t}_{\beta_{r+1}}$, is not zero, then the term $\diff{t}_{\beta_{r+1}}$
\emph{gives every contribution}, i.e. it is \emph{complete}.

The particular situation of the leading term $\diff{t}_{\beta_{1}}$
is more natural, and is tied to the idea that $\left(\diff{t}_{\beta_{1}}\right)^{\frac{1}{p}}=\diff{t}_{p\beta}$,
so that $\left(\diff{t}_{\beta_{1}}\right)^{\frac{1}{p}}=0$ if and
only if $\beta_{1}<\frac{1}{p}$.

All this motivates the following

\begin{defn}
\label{def:incompleteTerm} Under the hypotheses introduced at the beginning
of this section, we say that $\diff{t}_{\beta_{r}}$ \emph{loses 
information in} $x^{\frac{1}{p}}$, or that $\diff{t}_{\beta_{r}}$
\emph{is incomplete with respect to} $x^{\frac{1}{p}}$, if and only
if the following conditions hold:

\begin{enumerate}
\item \label{enu:1_defIncompleteTerm}$r=1\then\beta_{1}<\frac{1}{p}$
\item \label{enu:2_defIncompleteTerm}If $r>1$, then
\[
\exists n=1,\dots,k_{x,\frac{1}{p}}\ \exists\gamma\in\N^{N-1}:\ |\gamma|=n\ ,\ \gamma_{r-1}\ne0\ ,\ \omega(x,\gamma)\oplus p\beta_{1}<1.
\]

\end{enumerate}

\end{defn}

Let us analyze the condition \ref{enu:2_defIncompleteTerm} to see
that it corresponds to our intuition:

\begin{itemize}
\item '$\exists n=1,\dots,k_{x,\frac{1}{p}}\ \exists\gamma\in\N^{N-1}:\ |\gamma|=n$':
looking at \eqref{eq:x_1over_p_powerFormula}, we can say that this
part of \ref{enu:2_defIncompleteTerm} corresponds to {}``at least
one summand in the expansion of $x^{1/p}$''
\item '$\gamma_{r-1}\ne0$': {}``where the term $\diff{t}_{\beta_{r}}$
appears''
\item '$\omega(x,\gamma)\oplus p\beta_{1}<1$': {}``is zero''.

\end{itemize}

Consider the case $p=\frac{1}{q}$, where $q\in\N_{>0}$.
In this case, the power $x^{1/p}$ becomes $x^{q}$ and hence we 
obtain an equivalent formulation starting
from the multinomial formula:
\[
x^{q}=\left(\sum_{j=1}^{N}  b_{j}\diff{t}_{\beta_{j}}\right)^{q}=\sum_{\substack{\eta\in\N^{N}\\
|\eta|=q}
} \frac{q!}{\eta!} b_{1}^{\eta_{1}}\cdot\ldots\cdot b_{N}^{\eta_{N}}\diff{t}_{\bigoplus_{j=1}^{N}\frac{\beta_{j}}{\eta_{j}}}.
\]
In fact, we have

\begin{thm}
If $p=\frac{1}{q}$, $q\in\N_{>0}$, with $q\le\ulcorner\beta_{2}\ominus\beta_{1}\urcorner$,
then we have that
\[
\diff{t}_{\beta_{r+1}}\text{ loses information in }x^{q}
\]
if and only if

\begin{equation}
\exists\eta\in\N^{N}:\ |\eta|=q\ ,\ \eta_{r+1}\ne0\ ,\ \bigoplus_{j=1}^{N}\frac{\beta_{j}}{\eta_{j}}<1.\label{eq:equivalentUsingMultinomial}
\end{equation}

\end{thm}

\noindent \textbf{Proof:} We first compute the term $\omega(x,\gamma)\oplus\frac{\beta_{1}}{q}$
of Definition \ref{def:incompleteTerm}, in the case $|\gamma|=q-a$,
where $a\in\N_{\le q}$:

\begin{align*}
\left(\omega(x,\gamma)\oplus\frac{\beta_{1}}{q}\right)^{-1} & =\left(\frac{\beta_{1}}{q}\oplus\bigoplus_{j=2}^{N}\frac{\beta_{j}\ominus\beta_{1}}{\gamma_{j-1}}\right)^{-1}=\\
 & =\frac{q}{\beta_{1}}+\sum_{j=2}^{N}\frac{\gamma_{j-1}}{\beta_{j}\ominus\beta_{1}}=\\
 & =\frac{q}{\beta_{1}}+\sum_{j=2}^{N}\gamma_{j-1}\cdot\left(\frac{1}{\beta_{j}}-\frac{1}{\beta_{1}}\right)=\\
 & =\frac{q}{\beta_{1}}+\sum_{j=2}^{N}\frac{\gamma_{j-1}}{\beta_{j}}-\frac{q}{\beta_{1}}+\frac{a}{\beta_{1}},
\end{align*}
where the last equality is due to the hypothesis $|\gamma|=q-a$. Therefore

\begin{equation}
\omega(x,\gamma)\oplus\frac{\beta_{1}}{q}=\frac{\beta_{1}}{a}\oplus\bigoplus_{j=2}^{N}\frac{\beta_{j}}{\gamma_{j-1}}\quad\text{if }|\gamma|=q-a.\label{eq:1_omegaOplusBeta1OverQ}
\end{equation}

To prove that \eqref{eq:equivalentUsingMultinomial} is necessary,
we start from the hypothesis that there exist $n=1,\dots,k_{x,\frac{1}{p}}=\min\left(\ulcorner\beta_{2}\ominus\beta_{1}\urcorner,q\right)=q$
and there exists $\gamma\in\N^{N-1}$ such that
\[
|\gamma|=n\quad,\quad\gamma_{r}\ne0\quad,\quad\omega(x,\gamma)\oplus\frac{\beta_{1}}{q}<1.
\]
Set
\[
\eta_{j}:=
\begin{cases}
q-n & \text{ if }j=1\\
\gamma_{j-1} & \text{ if }j=2,\dots,N
\end{cases}
\]
Then $\eta\in\N^{N}$, $|\eta|=|\gamma|+q-n=q$, and $\eta_{r+1}=\gamma_{r}\ne0$.
Finally, from \eqref{eq:1_omegaOplusBeta1OverQ}, with $a:=q-n$,
we obtain
\[
\bigoplus_{j=1}^{N}\frac{\beta_{j}}{\eta_{j}}=\frac{\beta_{1}}{q-n}\oplus\bigoplus_{j=2}^{N}\frac{\beta_{j}}{\gamma_{j-1}}=\omega(x,\gamma)\oplus\frac{\beta_{1}}{q}<1,
\]
which concludes the first part of our proof.

To prove that \eqref{eq:equivalentUsingMultinomial} is sufficient
to obtain that $\diff{t}_{\beta_{r+1}}$ is incomplete, we consider
$\eta$ as in \eqref{eq:equivalentUsingMultinomial} and set $\gamma:=(\eta_{2},\dots,\eta_{N})\in\N^{N-1}$.
Then
\[
|\gamma|=|\eta|-\eta_{1}=q-\eta_{1}\le q=\min\left(\ulcorner\beta_{2}\ominus\beta_{1}\urcorner,q\right)=k_{x,\frac{1}{p}}
\]
and $\gamma_{r}=\eta_{r+1}\ne0$. Using \eqref{eq:1_omegaOplusBeta1OverQ}
with $a:=\eta_{1}\le|\eta|=q$, we have
\[
\omega(x,\gamma)\oplus\frac{\beta_{1}}{q}=\frac{\beta_{1}}{\eta_{1}}\oplus\bigoplus_{j=2}^{N}\frac{\beta_{j}}{\gamma_{j-1}}=\bigoplus_{j=1}^{N}\frac{\beta_{j}}{\eta_{j}}<1,
\]
which concludes our proof.$\qed$

The following theorem can be viewed as a validation of our definition of incomplete term.

\begin{thm}
\label{thm:rootsAreInvertibleOnCompleteNumbers}Let $x\in D_{\infty}$
be a strictly positive infinitesimal, and let $p\in\R_{>0}$. Suppose
that in the decomposition of $x$ no term loses information in
$x^{\frac{1}{p}}$, then
\[
\left(x^{\frac{1}{p}}\right)^{p}=x.
\]
Therefore, the power $(-)^{p}$ is an injection on the set
\[
\left\{ x\in D_{\infty}\,|\, \text{all terms of } x \text{ are }(-)^{\frac{1}{p}}\text{ complete}\right\} \cup\left\{ y\in\ER\,|\,\st{y}\ne0\right\} .
\]

\end{thm}

\noindent\textbf{Proof:} If $p>1$, setting $q:=\frac{1}{p}<1$, 
the conclusion follows from Theorem \ref{thm:inverseUpToInf}: $\left(x^{q}\right)^{\frac{1}{q}}=x=\left(x^{\frac{1}{p}}\right)^{p}.$ Therefore,
the only interesting case is $p<1$.

If, in the decomposition of $x$, we have only $N=1$ term, then,
by hypothesis, this term is not incomplete. Taking the negation of
Definition \ref{def:incompleteTerm} for $r=1$, we obtain $\beta_{1}\ge\frac{1}{p}$
and
\[
\left(x^{\frac{1}{p}}\right)^{p}=\left[\left(b_{1}\diff{t}_{\beta_{1}}\right)^{\frac{1}{p}}\right]^{p}=\left(b_{1}^{\frac{1}{p}}\diff{t}_{p\beta_{1}}\right)^{p}=b_{1}\diff{t}_{\beta_{1}}=x.
\]
Let us observe that since $p\beta_{1}\ge1$, the expression $b_{1}^{\frac{1}{p}}\diff{t}_{p\beta_{1}}$
is a decomposition, so that in taking its $p$-th power, we have applied
Definition \ref{def:c^p}. We can hence suppose $N>1$.

As usual, we refer to \eqref{eq:x_1over_p_powerFormula}. By hypothesis,
every term $\diff{t}_{\beta_{j}}$ is complete, which means
\[
p\beta_{1}\ge1
\]
and, for every $j=2,\dots,N$, $n=1,\dots,k$, and $\gamma\in\N^{N-1}$
we must have
\[
|\gamma|=n\ ,\ \gamma_{j-1}\ne0\then\omega(x,\gamma)\oplus p\beta_{1}\ge1.
\]
This implies that in \eqref{eq:x_1over_p_powerFormula} the greatest
infinitesimal is $\diff{t}_{p\beta_{1}}$ and every summand is not
zero. Therefore, to compute $\left(x^{\frac{1}{p}}\right)^{p}$, we
can use Lemma \ref{lem:formulaForPowersWithoutUsingDecomp}:

\begin{align*}
\left(x^{\frac{1}{p}}\right)^{p} & =b_{1}\diff{t}_{\beta_{1}}\cdot\left(1+\sum_{n=1}^{k}\binom{\frac{1}{p}}{n}\sum_{|\gamma|=n}\frac{n!}{\gamma!}\cdot b_{1}^{-n}\cdot b_{2}^{\gamma_{1}}\cdot\ldots\cdot b_{N}^{\gamma_{N-1}}\diff{t}_{\omega(x,\gamma)}\right)^{p}\\
 & =:b_{1}\diff{t}_{\beta_{1}}\cdot\left(1+h\right)^{p}.
\end{align*}

We have to prove that

\begin{align*}
b_{1}\diff{t}_{\beta_{1}}\cdot\left(1+h\right)^{p} & =x\\
 & =b_{1}\diff{t}_{\beta_{1}}+b_{2}\diff{t}_{\beta_{2}}+\ldots+b_{N}\diff{t}_{\beta_{N}}\\
 & =:b_{1}\diff{t}_{\beta_{1}}+k.
\end{align*}

To this end, we use the following

\begin{lem}
\label{lem:ratioOfInfinitesimals}Let $h$, $k\in D_{\infty}$ and
$p$, $b$, $\beta\in\R_{\ne0}$ such that $p<1$ and $\beta>\omega(k)$.
Define
\[
\frac{k}{b\diff{t}_{\beta}}:=\sum_{i=1}^{N}\frac{a_{i}}{b}\diff{t}_{\alpha_{i}\ominus\beta},
\]
where $k=\sum_{i=1}^{N}a_{i}\diff{t}_{\alpha_{i}}$ is the decomposition
of $k$. Then
\begin{equation}
1+h=\left(1+\frac{k}{b\diff{t}_{\beta}}\right)^{\frac{1}{p}}\then b\diff{t}_{\beta}\cdot(1+h)^{p}=b\diff{t}_{\beta}+k.\label{eq:implicationRatioOfInfinitesimals}
\end{equation}

\end{lem}

\noindent Therefore, to obtain our conclusion it suffices to verify the assumption of \eqref{eq:implicationRatioOfInfinitesimals},
that is
\begin{align*}
1+\sum_{n=1}^{k}\binom{\frac{1}{p}}{n}\sum_{|\gamma|=n}\frac{n!}{\gamma!}\cdot b_{1}^{-n}\cdot b_{2}^{\gamma_{1}}\cdot\ldots\cdot b_{N}^{\gamma_{N-1}}\diff{t}_{\omega(x,\gamma)} & =\\
=\left(1+\frac{b_{2}}{b_{1}}\diff{t}_{\beta_{2}\ominus\beta_{1}}+\ldots+\frac{b_{N}}{b_{1}}\diff{t}_{\beta_{N}\ominus\beta_{1}}\right)^{\frac{1}{p}}.
\end{align*}
Indeed,
\noindent 
\begin{multline*}
\left(1+\frac{b_{2}}{b_{1}}\diff{t}_{\beta_{2}\ominus\beta_{1}}+\ldots+\frac{b_{N}}{b_{1}}\diff{t}_{\beta_{N}\ominus\beta_{1}}\right)^{\frac{1}{p}}=1+\sum_{n=1}^{k}\binom{\frac{1}{p}}{n}\left(\sum_{j=2}^{N}\frac{b_{j}}{b_{1}}\diff{t}_{\beta_{j}\ominus\beta_{1}}\right)^{n}\\
=1+\sum_{n=1}^{k}\binom{\frac{1}{p}}{n}\sum_{|\gamma|=n}\frac{n!}{\gamma!}\cdot\left(\frac{b_{2}}{b_{1}}\right)^{\gamma_{1}}\cdot\ldots\cdot\left(\frac{b_{N}}{b_{1}}\right)^{\gamma_{N-1}}\diff{t}_{\bigoplus_{j=2}^{N}\frac{\beta_{j}\ominus\beta_{1}}{\gamma_{j-1}}}\\
=1+\sum_{n=1}^{k}\binom{\frac{1}{p}}{n}\sum_{|\gamma|=n}\frac{n!}{\gamma!}\cdot b_{1}^{-n}\cdot b_{2}^{\gamma_{1}}\cdot\ldots\cdot b_{N}^{\gamma_{N-1}}\diff{t}_{\omega(x,\gamma)}=1+h,
\end{multline*}
as claimed.
$\qed$

\noindent \textbf{Proof of Lemma \ref{lem:ratioOfInfinitesimals}:}
It suffices to note that the powers $(-)^{p}$ and $(-)^{\frac{1}{p}}$
are smooth if applied to invertible Fermat reals. By the elementary
transfer theorem, they hence have all the usual properties, so that
we can write
\[
\left(1+h\right)^{p}=\left[\left(1+\frac{k}{b\diff{t}_{\beta}}\right)^{\frac{1}{p}}\right]^{p}=1+\frac{k}{b\diff{t}_{\beta}}.
\]
Thus,
\begin{align*}
b\diff{t}_{\beta}\cdot\left(1+h\right)^{p} & =b\diff{t}_{\beta}+b\diff{t}_{\beta}\cdot\frac{k}{b\diff{t}_{\beta}}=\\
 & =b\diff{t}_{\beta}+b\diff{t}_{\beta}\cdot\left(\sum_{i=1}^{N}\frac{a_{i}}{b}\diff{t}_{\alpha_{1}\ominus\beta}\right)=\\
 & =b\diff{t}_{\beta}+k 
\end{align*}
$\qed$

\begin{rem}
To simplify the notations, let us define
\[
x\text{ is }(-)^{\frac{1}{p}}\text{-complete}\quad:\iff\quad\forall j=1,\dots,N:\ \diff{t}_{\beta_{j}}\text{ is complete w.r.t. }x^{\frac{1}{p}}
\]
\[
\text{Compl}\left(\frac{1}{p}\right):=\left\{ x\in D_{\infty}\,|\, x\text{ is }(-)^{\frac{1}{p}}\text{-complete}\right\} \cup\left\{ y\in\ER\,|\,\st{y}\ne0\right\} .
\]
Then the map
\[
(-)^{p}|_{\text{Compl}\left(\frac{1}{p}\right)}:\text{Compl}\left(\frac{1}{p}\right)\freccia\ER
\]
is injective. Moreover, the map
\[
(-)^{p}:\ER\freccia\ER
\]
is surjective. However, the power $(-)^{p}$ doesn't map $\text{Compl}\left(\frac{1}{p}\right)$
onto itself. In fact, if $p=\frac{1}{2}$ and $y=\diff{t}_{2}$, then
$y^{\frac{1}{p}}=y^{2}=\diff{t}$ so that every summand in $y^{\frac{1}{p}}$
is nonzero and hence $y\in\text{Compl}\left(\frac{1}{p}\right)$.
This $y$ is therefore in the codomain, but it is not of the form
$y=x^{p}$ for $x\in\text{Compl}\left(\frac{1}{p}\right)$,
as otherwise $y^{2}=\left(x^{\frac{1}{2}}\right)^{2}=x$,
so that $x=\diff{t}$ and we would also have $\diff{t}\in\text{Compl}\left(2\right)$,
contradicting $\left(\diff{t}\right)^{2}=0$.

\end{rem}

We can now prove that complete Fermat infinitesimals have very favorable
properties related to powers.

\begin{thm}
\label{thm:completeAndProp_xTo_pq}Let $x\in D_{\infty}$ be a strictly
positive infinitesimal, and $p$, $q\in\R_{>0}$. Then
\[
x\text{ is }(-)^{p}\text{-complete}\then\left(x^{p}\right)^{q}=x^{pq}.
\]

\end{thm}
\noindent \textbf{Proof:} Let us start from the usual formula \eqref{eq:finalFormula}
applied to $x^{p}$:
\[
x^{p}=b_{1}^{p}\diff{t}_{\frac{\beta_{1}}{p}}+\sum_{n=1}^{k_{x,p}}\binom{p}{n}\sum_{\gamma\in\Gamma_{x,n}}\frac{n!}{\gamma!}\cdot b_{1}^{p-n}\cdot b_{2}^{\gamma_{1}}\cdot\ldots\cdot b_{N}^{\gamma_{N-1}}\diff{t}_{\omega(x,\gamma)\oplus\frac{\beta_{1}}{p}}.
\]
By hypothesis, $x$ is $(-)^{\frac{1}{p}}$-complete, that is $\frac{\beta_{1}}{p}\ge1$,
and for every $j=2,\dots,N$, $n=1,\dots,k_{x,p}$, and $\gamma\in\N^{N-1}$
we have

\begin{equation}
|\gamma|=n\ ,\ \gamma_{j-1}\ne0\then\omega(x,\gamma)\oplus\frac{\beta_{1}}{p}\ge1.\label{eq:everySummandIsNotZero}
\end{equation}

This property implies that $\Gamma_{x,n}$ is trivial, in fact, in
general
\[
\forall a,b\in\R:\ a<1\then\left(a\oplus b\right)^{-1}=\frac{1}{a}+\frac{1}{b}>\frac{1}{a}>1.
\]
Therefore
\[
\forall n=1,\dots,k_{x,p}\ \forall\gamma\in\N^{N-1}:\ |\gamma|=n\then\omega(x,\gamma)\ge1,
\]

\begin{equation}
x^{p}=b_{1}^{p}\diff{t}_{\frac{\beta_{1}}{p}}+\sum_{n=1}^{k_{x,p}}\binom{p}{n}\cdot\sum_{\substack{\gamma\in\N^{N-1}\\
|\gamma|=n}
}\frac{n!}{\gamma!}\cdot b_{1}^{p-n}\cdot b_{2}^{\gamma_{1}}\cdot\ldots\cdot b_{N}^{\gamma_{N-1}}\diff{t}_{\omega(x,\gamma)\oplus\frac{\beta_{1}}{p}}.\label{eq:xpWithoutGamma}
\end{equation}

By \eqref{eq:everySummandIsNotZero}, every summand in \eqref{eq:xpWithoutGamma}
has order greater or equal to 1. We can hence apply Lemma \ref{lem:formulaForPowersWithoutUsingDecomp}
obtaining

\begin{align*}
\left(x^{p}\right)^{q} & =b_{1}^{pq}\diff{t}_{\frac{\beta_{1}}{pq}}\cdot\left(1+\sum_{n=1}^{k_{x,p}}\binom{p}{n}\sum_{\substack{\gamma\in\N^{N-1}\\
|\gamma|=n}
}\frac{n!}{\gamma!}\cdot b_{1}^{-n}\cdot b_{2}^{\gamma_{1}}\cdot\ldots\cdot b_{N}^{\gamma_{N-1}}\diff{t}_{\omega(x,\gamma)}\right)^{q}\\
 & =b_{1}^{pq}\diff{t}_{\frac{\beta_{1}}{pq}}\cdot\left[\left(1+\sum_{j=2}^{N}\frac{b_{j}}{b_{1}}\diff{t}_{\beta_{j}\ominus\beta_{1}}\right)^{p}\right]^{q}\\
 & =b_{1}^{pq}\diff{t}_{\frac{\beta_{1}}{pq}}\cdot\left(1+\sum_{j=2}^{N}\frac{b_{j}}{b_{1}}\diff{t}_{\beta_{j}\ominus\beta_{1}}\right)^{pq}\\
 & =x^{pq}.
\end{align*}
Here the property $\left[(1+h)^{p}\right]^{q}=\left(1+h\right)^{pq}$,
for $h\in D_{\infty}$, holds because the base is invertible.$\qed$

The following result does not depend on the notion of complete term,
but supposes that in the product $x\cdot y$ no term becomes zero.

\begin{thm}
Let $x$, $y\in\ER_{>0}$ and $p\in\R_{>0}$, such that

\begin{equation}
\forall i=1,\dots,N_{y}\ \forall j=1,\dots,N_{x}:\ \omega_{i}(y)\oplus\omega_{j}(x)\ge1,\label{eq:noTermsZeroInProduct}
\end{equation}
where $N_{x}$ and $N_{y}$ are the number of summands in the decompositions
of $x$, $y$ respectively. Then
\[
\left(x\cdot y\right)^{p}=x^{p}\cdot y^{p}.
\]

\end{thm}

Let us observe that the hypothesis \eqref{eq:noTermsZeroInProduct}
reduces to the usual $(-)^{2}$-com\-ple\-te\-ness in case $x=y$,
i.e. in case of the property $\left(x^{2}\right)^{p}=x^{p}\cdot x^{p}$.

\noindent \textbf{Proof:} We will proceed for $x$, $y\in D_{\infty}$,
the proof being analogous if $x$ or $y$ is invertible. Let $x=\sum_{j=1}^{N}b_{j}\diff{t}_{\beta_{j}}$
and $y=\sum_{i=1}^{M}a_{i}\diff{t}_{\alpha_{i}}$ be the decompositions
of $x$ and $y$. Then

\begin{equation}
x\cdot y=\sum_{i,j}b_{j}a_{i}\diff{t}_{\alpha_{i}\oplus\beta_{j}}.\label{eq:product}
\end{equation}
By hypothesis, every summand in this sum is nonzero. Of course, $\alpha_{1}\oplus\beta_{1}$
is the leading term in \eqref{eq:product}, and we can hence apply
Lemma \ref{lem:formulaForPowersWithoutUsingDecomp} to obtain
\begin{equation}
\left(xy\right)^{p}=a_{1}^{p}b_{1}^{p}\diff{t}_{\frac{\alpha_{1}\oplus\beta_{1}}{p}}\cdot\left(1+\sum_{\substack{i,j\\
(i,j)\ne(1,1)}
}\frac{b_{j}a_{i}}{b_{1}a_{1}}\diff{t}_{\alpha_{i}\oplus\beta_{j}\ominus\left(\alpha_{1}\oplus\beta_{1}\right)}\right)^{p}.\label{eq:xp_To_p}
\end{equation}
On the other hand
\begin{align}
x^{p} & =b_{1}^{p}\diff{t}_{\frac{\beta_{1}}{p}}\cdot\left(1+\sum_{j=2}^{N}\frac{b_{j}}{b_{1}}\diff{t}_{\beta_{j}\ominus\beta_{1}}\right)^{p}\label{eq:x_To_p}\\
y^{p} & =a_{1}^{p}\diff{t}_{\frac{\alpha_{1}}{p}}\cdot\left(1+\sum_{i=2}^{M}\frac{a_{i}}{a_{1}}\diff{t}_{\alpha_{j}\ominus\alpha_{1}}\right)^{p}.\label{eq:y_To_p}
\end{align}

The equality \eqref{eq:xp_To_p} can also be written as

\begin{align*}
\left(xy\right)^{p}=a_{1}^{p}b_{1}^{p}\diff{t}_{\frac{\alpha_{1}\oplus\beta_{1}}{p}}\cdot & \left(1+\sum_{j}\frac{b_{j}a_{1}}{b_{1}a_{1}}\diff{t}_{\alpha_{1}\oplus\beta_{j}\ominus\left(\alpha_{1}\oplus\beta_{1}\right)}\right.\\
 & \h{0.215}+\sum_{i}\frac{b_{1}a_{i}}{b_{1}a_{1}}\diff{t}_{\alpha_{i}\oplus\beta_{1}\ominus\left(\alpha_{1}\oplus\beta_{1}\right)}\\
 & \h{0.215}\left.+\sum_{\substack{i\ge2\\
j\ge2}
}\frac{b_{j}a_{i}}{b_{1}a_{1}}\diff{t}_{\alpha_{i}\oplus\beta_{j}\ominus\left(\alpha_{1}\oplus\beta_{1}\right)}\right)^{p},
\end{align*}
which is the same result obtained from multiplying \eqref{eq:x_To_p} and
\eqref{eq:y_To_p}.$\qed$

\subsection{Roots are not smooth}

Here we prove that the power function
\[
(-)^{p}:\ER_{>0}\freccia\ER\quad,\quad p<1
\]
is not (non standard) smooth on $D_{\infty}\setminus\{0\}$ (for the
notion of non standard smoothness, see \cite{Gio10e,Gio09}). This
is a naturally expected result, because the corresponding derivative
should be
\[
\frac{\diff{x^{p}}}{\diff{x}}(h)=p\cdot h^{p-1}=\frac{p}{h^{1-p}},
\]
and, intuitively, $\frac{1}{h^{1-p}}$ is an infinite, whereas in
$\ER$ we obviously do not have infinities. Therefore, the theory
of Fermat reals should be sufficiently complete to prove that the
power function $(-)^{p}$ is not smooth at $h\in D_{\infty}\setminus\{0\}$.
Indeed, from the generalized Taylor formula (see Theorem 12.1.3
and Definition 12.2.7 in \cite{Gio09}), if we assume that $(-)^{p}$
is smooth, we would have
\begin{equation}
\exists m\in\ER^{\N}:\ \forall k\in D_{\infty}:\ (h+k)^{p}=\sum_{j=0}^{+\infty}\frac{h^{j}}{j!}\cdot m_{j}.\label{eq:genTaylorFormula}
\end{equation}
For the sake of completeness, we recall that this sequence $m=\left(m_{j}\right)_{j\in\N}$
is unique up to first order infinitesimals, i.e. if $\bar{m}=\left(\bar{m_{j}}\right)_{j\in\N}$
verifies \eqref{eq:genTaylorFormula}, then $m_{j}=_{1}\bar{m}_{j}$
for every $j\in\N$.

\noindent Set $s:=-\text{sgn}(h)\in\{+1,-1\}$, and take $k=s\diff{t}$
in \eqref{eq:genTaylorFormula}. We have that $h+s\diff{t}\ne0$ because
$\text{sgn}(h)\ne\text{sgn}(s\diff{t})$, and taking into account that $m_{0}=0^{p}=0$,
we obtain
\[
\left(h+s\diff{t}\right)^{p}=m_{0}+s\diff{t}\cdot m_{1}=sm_{1}\diff{t}\in I_{1}=\{x\in\ER\,|\,\omega(x)\le1\}.
\]
Now, the order of $\left(h+s\diff{t}\right)^{p}$ is $\frac{\max(\omega(h),1)}{p}=\frac{\omega(h)}{p}$
because $h+s\diff{t}\ne0$ and since always $\omega(h)\ge1$. Therefore,
we would have $\frac{\omega(h)}{p}\le1$, i.e. $p\ge\omega(h)\ge1$,
whereas we have supposed $p<1$.

\subsection{An application to the infinitesimal Taylor formula with fractional
derivatives}

Using powers $h^{p}$ of infinitesimals $h\in D_{\infty}$, we can
prove an infinitesimal Taylor formula with fractional derivatives in a straightforward manner.
This further underlines the ease of translating classical results using
the infinitesimal language of the ring of Fermat reals. Frequently,
these translations are really faithful to the informal use sometimes
appearing in applications. Let us note that the same translations
are not so easily performed in algebraic models of infinitesimals,
like in Synthetic Differential Geometry (see, e.g., \cite{Mo-Re},
\cite{Koc}) or in Levi-Civita fields (\cite{Sh-Be,Sha})
or Weil functors (\cite{Ko-Mi-Sl,Kr-Mi}).

We start with some definitions and a theorem, taken from \cite{Od-Sh}.

\begin{defn}
\label{def:C_alpha_space}If $\alpha\in\R$, we will denote with $\mathcal{C}_{\alpha}(\R_{>0},\R)$
the set of all the functions $f:\R_{>0}\freccia\R$ that can be written
as
\[
f(x)=x^{p}\cdot f_{1}(x)\quad\forall x\in\R_{>0},
\]
for some $p>\alpha$ and some continuous function $f_{1}\in\mathcal{C}^{0}(\R_{>0},\R)$.
Moreover, for every $m\in\N_{>0}\cup\{\infty\}$ we also set
\[
\mathcal{C}_{\alpha}^{m}(\R_{>0},\R):=\left\{ f\in\R_{>0}\freccia\R\,|\, f^{(m)}\in\mathcal{C}_{\alpha}(\R_{>0},\R)\right\} .
\]

\end{defn}

Secondly, we define the Riemann-Liouville integral operator of order
$\alpha>0$ with $a\ge0$.

\begin{defn}
\label{def:RLOperator}Let $\alpha\in\R_{>0}$, $a\in\R_{\ge0}$ and
$f\in\mathcal{C}_{\alpha}(\R_{>0},\R)$, then

\begin{align}
J_{a}^{0}f(x) & :=f(x)\nonumber \\
J_{a}^{\alpha}f(x) & :=\frac{1}{\Gamma(\alpha)}\int_{a}^{x}(x-t)^{\alpha-1}f(t)\diff{t}\quad\forall x\in\R_{>a}.\label{eq:RLOperatorAlpha}\end{align}

\end{defn}
Here $\Gamma$ denotes the gamma function. To derive the fractional
Taylor formula, we need the Caputo fractional derivative.
\begin{defn}
\label{def:CaputoFractionalDerivative}Let $\alpha\in\R_{>0}$, $a\in\R_{\ge0}$,
and $f\in\mathcal{C}_{-1}^{m}(\R_{>0},\R)$. For simplicity of notations,
let $m:=\ulcorner\alpha\urcorner$ be the ceiling of $\alpha$. Then
\begin{align}
D_{a}^{\alpha}f & :\R_{\ge a}\freccia\R\label{eq:Caputo}\\
D_{a}^{\alpha}f(x) & :=J_{a}^{m-\alpha}f^{(m)}(x)\quad\forall x\ge a.
\end{align}
Finally, we set
\[
D_{a}^{n,\alpha}=D_{a}^{\alpha}\circ\ptind^{n}\circ D_{a}^{\alpha}\quad\forall n\in\N_{>0}.
\]

\end{defn}

\noindent From \eqref{eq:Caputo} and \eqref{eq:RLOperatorAlpha}
we therefore have
\[
D_{a}^{\alpha}f(x)=\frac{1}{\Gamma(m-\alpha)}\int_{a}^{x}(x-t)^{m-\alpha-1}f^{(m)}(t)\diff{t},
\]
where $m=\ulcorner\alpha\urcorner$.

The non-infinitesimal version of the generalized Taylor formula
with fractional derivatives is the following. For its proof, see \cite{Od-Sh}.

\begin{thm}
\label{thm:GTF_fractionalDerivatives}Let $\alpha$, $a$, $b\in\R$,
and $n\in\N$, with $0\le a<b$ and $0<\alpha\le1$. Consider a continuous
function $f\in\mathcal{C}_{0}([a,b],\R)$ such that
\[
D_{a}^{k,\alpha}f\in\mathcal{C}_{0}([a,b],\R)\quad\forall k=0,\dots,n+1.
\]
Then for every $x\in(a,b]$ there exists $\xi\in[a,x]$ such that
\[
f(x)=\sum_{i=0}^{n}\frac{(x-a)^{i\alpha}}{\Gamma(i\alpha+1)}D_{a}^{i,\alpha}f(a)+\frac{D_{a}^{n+1,\alpha}f(\xi)}{\Gamma\left((n+1)\alpha+1\right)}\cdot(x-a)^{(n+1)\alpha}.
\]

\end{thm}

In our framework we are able to prove a corresponding infinitesimal
Taylor formula for the following class of smooth functions:

\begin{defn}
\label{def:smoothForGTF_fractional}In the hypothesis of the previous
theorem, we set
\[
\mathcal{C}_{\alpha}^{\infty}([a,b],\R):=\left\{ f\in\mathcal{C}^{\infty}([a,b],\R)\,|\, D_{a}^{k,\alpha}f\in\mathcal{C}_{0}([a,b],\R)\quad\forall k\in\N\right\} .
\]

\end{defn}

\noindent Finally, we can state the main result of this section:

\begin{thm}
\label{thm:infTaylorFormulaFractional}Let $\alpha$, $a'$, $a$,
$b\in\R$, and $n\in\N$, with $0\le a'<a<b$ and $0<\alpha\le1$.
Consider a smooth function $f\in\mathcal{C}_{\alpha}^{\infty}([a,b],\R)$,
then
\[
\forall h\in D_{(n+1)\alpha - 1}:\ f(a+h)=\sum_{i=0}^{n}\frac{h^{i\alpha}}{\Gamma(i\alpha+1)}D_{a}^{i,\alpha}f(a).
\]

\end{thm}

\noindent \textbf{Proof:} Let $h=\sum_{j=1}^{N}b_{j}\diff{t}_{\beta_{j}}$
be the decomposition of the infinitesimal $h\in D_{(n+1)\alpha}$.
By Definition \ref{def:c^p}, we have
\[
h^{i\alpha}=b_{1}^{i\alpha}\diff{t}_{\frac{\beta_{1}}{i\alpha}}\cdot\left(1+\sum_{j=2}^{N}\frac{b_{j}}{b_{1}}\diff{t}_{\beta_{j}\ominus\beta_{1}}\right)^{i\alpha}\quad\forall i=0,\dots,n.
\]
This means, using an innocuous abuse of language, that $h_{t}:=\sum_{j=1}^{N}b_{j}t^{\frac{1}{\beta_{j}}}$
and
\[
h_{t}^{i\alpha}:=b_{1}^{i\alpha}t^{\frac{i\alpha}{\beta_{1}}}\cdot\left(1+\sum_{j=2}^{N}\frac{b_{j}}{b_{1}}t^{\frac{1}{\beta_{j}\ominus\beta_{1}}}\right)^{i\alpha}\quad\forall t\ge0
\]
are little-oh polynomials representing the Fermat real $h$ and $h^{i\alpha}$
respectively, and
\[
\left(h_{t}\right)^{i\alpha}=h_{t}^{i\alpha}\quad\forall t\ge0.
\]
For $t\ge0$ sufficiently small, we have $a'<a+h_{t}<b$, and we can
apply Theorem \ref{thm:GTF_fractionalDerivatives} at the point $a$,
obtaining

\begin{equation}
f(a+h_{t})=\sum_{i=0}^{n}\frac{h_{t}^{i\alpha}}{\Gamma(i\alpha+1)}D_{a}^{i,\alpha}f(x)+\frac{D_{a}^{n+1,\alpha}f(\xi_{t})}{\Gamma\left((n+1)\alpha+1\right)}\cdot h_{t}^{(n+1)\alpha}\quad\xi_{t}\in[a,a+h_{t}].\label{eq:TaylorWithXi}
\end{equation}

Now, $h\in D_{(n+1)\alpha - 1}$ so that $\omega(h)<(n+1)\alpha$ and
$h^{(n+1)\alpha}=0$, that is

\begin{equation}
\lim_{t\to0^{+}}\frac{h_{t}^{(n+1)\alpha}}{t}=0.\label{eq:hPowerIsZero}
\end{equation}

Considering that $D_{a}^{n+1,\alpha}f$ is continuous on $[a,b]$
and that $\xi_{t}\in[a,b]$, from \eqref{eq:hPowerIsZero} and \eqref{eq:TaylorWithXi}
we obtain
\[
f(a+h_{t})=\sum_{i=0}^{n}\frac{h_{t}^{i\alpha}}{\Gamma(i\alpha+1)}D_{a}^{i,\alpha}f(x)+o(t)\quad\text{as }t\to0^{+},
\]
yielding the claim.$\qed$

\section{Computer implementation}

The definition of the ring of Fermat reals is highly constructive.
Therefore, using object oriented programming, it is not hard to write
a computer code corresponding to $\ER$. We realized a first version
of this software using Matlab R2010b.

The constructor of a Fermat real is \texttt{x=FermatReal(s,w,r)},
where \texttt{s} is the $n+1$ double vector of standard parts (\texttt{s(1)}
is the standard part $\st{x}$) and \texttt{w} is the double vector
of orders (\texttt{w(1)} is the order $\omega(x)$ if $x\in\ER\setminus\R$,
otherwise \texttt{w={[}{]}} is the empty vector). The last input \texttt{r}
is a logical variable and assumes the value \texttt{true} if we want that
the display of the number \texttt{x} is realized using the Matlab
\texttt{rats} function for both its standard parts and orders. In
this way, the number will be displayed using continued fraction approximations
and therefore, in many cases, the calculations will be exact. These
inputs are the basic methods of every Fermat real, and can be accessed
using the \texttt{subsref}, and \texttt{subsasgn}, notations \texttt{x.stdParts},
\texttt{x.orders}, \texttt{x.rats}. The function \texttt{w=orders(x)}
gives exactly the double vector \texttt{x.orders} if $x\in\ER\setminus\R$
and \texttt{0} otherwise.

The function \texttt{dt(a)}, where \texttt{a} is a double, constructs
the Fermat real $\diff{t}_{a}$. Because we have overloaded all the
algebraic operations, like \texttt{x+y}, \texttt{x{*}y}, \texttt{x-y},
\texttt{-x}, \texttt{x==y}, \texttt{x\textasciitilde{}=y}, \texttt{x<y},
\texttt{x<=y}, \texttt{x\textasciicircum{}y}, we can define a Fermat
real e.g. using an expression of the form \texttt{x=2+3{*}dt(2)-1/3{*}dt(1)},
which corresponds to \texttt{x=FermatReal({[}2 3 -1/3{]},{[}2 1{]},true)}.

We have also realized the function \texttt{y=decomposition(x)}, which
gives the decomposition of the Fermat real \texttt{x}, \texttt{abs(x)},
\texttt{log(x)}, \texttt{exp(x)}, \texttt{isreal(x)}, \texttt{is\-in\-fi\-ni\-te\-si\-mal(x)},
\texttt{isinvertible(x)}.

The logical function \texttt{v=eqUpTo(k,x,y)} corresponds to $x=_{k}y$.

The ratio \texttt{x/y} (see Theorem \ref{thm:linearEquations}) has
been implemented for \texttt{x} and \texttt{y} infinitesimals and
\texttt{y\textasciitilde{}=0}, or in case \texttt{y} is invertible.
Finally, the function \texttt{y=ext(f,x)}, corresponds to $\ext{f}(x)$
and has been realized using the evaluation of the symbolic Taylor
formula of the inline function \texttt{f}.

The functions \texttt{dF} and \texttt{dOmega} correspond, respectively,
to the Fermat and the omega distance, while \texttt{x\textasciicircum{}p},
\texttt{sqrt(x)} and \texttt{nthroot(x,n)} have been realized both
for \texttt{x} infinitesimal or invertible using the formulas we have
derived in the present work.

Using these tools, we can easily find, e.g., that
\[
\frac{\sin(\sqrt{\diff{t}_{3}+2\diff{t}_{2}})}{\cos(\sqrt[3]{-4\diff{t}})}=\diff{t}_{6}+\diff{t}_{3}-\frac{2}{3}\diff{t}_{2}+\frac{1096}{2787}\diff{t}_{\frac{6}{5}}+\frac{1234}{913}\diff{t}.
\]
This corresponds to the following Matlab code:

\noindent \texttt{>\textcompwordmark{}> x=sqrt(dt(3)+2{*}dt(2))}

\texttt{x = }

\texttt{dt\_6 + dt\_3 - 1/2{*}dt\_2 + 1/2{*}dt\_3/2 - 5/8{*}dt\_6/5}

\noindent \texttt{>\textcompwordmark{}> y=nthroot(-4{*}dt(1),3)}

\texttt{y = }

\texttt{-1008/635{*}dt\_3}

\noindent \texttt{>\textcompwordmark{}> g=inline('cos(y)')}

\texttt{g =}

\texttt{Inline function: g(y) = cos(y)}

\texttt{>\textcompwordmark{}> f=inline('sin(x)')}

\texttt{f =}

\texttt{Inline function: f(x) = sin(x)}

\noindent \texttt{>\textcompwordmark{}> decomposition(ext(f,x)/ext(g,y))}

\texttt{ans = }

\texttt{dt\_6 + dt\_3 - 2/3{*}dt\_2 + 1096/2787{*}dt\_6/5 + 1234/913{*}dt}

Up to now, this code has been written only to show concretely the possibilities of the ring $\ER$. On the other hand, it is clear that it is possible to write it with a more specific aim. For example, as in case of the Levi-Civita field (\cite{Berz-et-al,Sha}) possible applications of a specifically rewritten code include automatic differentiation theory. Let us note that, even if the theory of Fermat reals applies to smooth functions, a full treatment of right and left sided derivatives is possible (\cite{Gio09}), so that the theory can be applied consistently also to piecewise smooth functions. Finally, the use of nilpotent elements permits to fully justify that every derivative estimation of a computer function (\cite{Sha}) reduces to a finite number of algebraic calculations.

The Matlab source code is freely available under open-source licence,
and can be requested from the authors of the present article.

\section{Conclusions}

Usually, it is common to study extended structures, like the ring
of Fermat reals $\ER$, using suitable extensions of well established
notions. For example, it is more natural to search for metrics of
the form $d:\ER\times\ER\freccia\ER$ than for standard metrics on the
set $\ER$. We have shown that it is possible, and also very natural,
to define standard topological structures on the ring $\ER$ having
very favorable relationships with various aspects of differential calculus
on $\ER$. This allows a better dialog with mathematicians not already
familiar with the theory of Fermat reals and underlines that this ring
is not {}``non standard''.

Moreover, with the present work, we are continuing our program to
define a meaningful and powerful ring with infinitesimals
using only very well behaved representative functions for new numbers.
If one thinks at non standard analysis or Colombeau's ring of generalized
numbers, it becomes clear that this quest is nontrivial.
As a consequence, we have been able to characterize ideals of the ring $\ER$ 
in a very simple and descriptive way.

Finally, we have proved that nilpotent elements and arbitrary roots
can coexist very well, even if this seems impossible at a first glance.
This is a very important step toward the idea of using nilpotent infinitesimals
for stochastic calculus. For example, based on a very helpful discussion
with N. Blagowest (Department of Physics, K. Preslawki University,
Bulgaria) we may call \emph{Ito process} any (deterministic) function
$x:\ER\freccia\ER$ such that
\[
\forall t\in\R\ \forall h\in D_{\frac{1}{2}}:\ x(t+h)=x(t)+v\left[t,x(t)\right]\cdot h+\lambda\sqrt{h}.
\]
In this approach, the deep mathematical problem is that there doesn't
exist a non trivial smooth function that verifies such a definition.
Of course, we need continuous but nowhere differentiable functions
and hence, we need to extend the ring $\ER$ by suitable infinities.
Indeed, this is indispensable if we want to consider the derivatives
of the function $x$. Therefore, the new problem to face becomes:
can infinities coexist as reciprocals of nilpotent infinitesimals?
This question will be the subject of future work.

\clearpage{}

\end{document}